\documentclass[aap,MSNbibl,citesort,dvips]{arximspdf}
\usepackage{dcolumn}
\usepackage{shere}
\usepackage{graphicx}

\doi{10.1214/10-AAP711}
\volume{21}
\issue{3}
\pubyear{2011}
\firstpage{908}
\lastpage{950}

\makeatletter

\newproclaim{assumption}{Assumption}
\newtheorem{theorem}{Theorem}
\newtheorem{lemma}{Lemma}

\newcolumntype{d}[1]{D{.}{.}{#1}}

\makeatother

\begin{document}
\begin{frontmatter}

\title{Efficient simulation of nonlinear parabolic SPDEs with additive noise}
\runtitle{Efficient simulation of SPDEs}

\begin{aug}
\author[A]{\fnms{Arnulf} \snm{Jentzen}\thanksref{T0,T1}\ead[label=e1]{ajentzen@math.princeton.edu}},
\author[B]{\fnms{Peter} \snm{Kloeden}\corref{}\thanksref{T1}\ead[label=e2]{kloeden@math.uni-frankfurt.de}} and
\author[B]{\fnms{Georg} \snm{Winkel}\ead[label=e3]{georg.winkel@gmx.net}}
\runauthor{A. Jentzen, P. Kloeden and G. Winkel}
\affiliation{Princeton University, Johann Wolfgang Goethe
University
and~Johann~Wolfgang~Goethe University}
\address[A]{A. Jentzen\\
Program in Applied\\
\quad and Computational Mathematics\\
Princeton University\\
Fine Hall, Washington Road\\
Princeton, New Jersey 08544-1000\\
USA\\
\printead{e1}}
\address[B]{P. Kloeden\\
G. Winkel\\
Institute of Mathematics\\
Johann Wolfgang Goethe University\\
D-60054 Frankfurt am Main\\
Germany\\
\printead{e2}\\
\phantom{E-mail: }\printead*{e3}}
\end{aug}

\thankstext{T0}{Supported in part by the DFG Collaborative Research Centre 701
``Spectral Structures and Topological Methods in\vspace*{1pt} Mathematics.''}

\thankstext{T1}{Supported by the DFG project ``Pathwise numerics and
dynamics of stochastic evolution equations.''}

\received{\smonth{11} \syear{2009}}
\revised{\smonth{4} \syear{2010}}

%
\begin{abstract}
Recently, in
a paper by Jentzen and Kloeden
 [\textit{Proc. R. Soc. Lond. Ser.~A Math. Phys. Eng. Sci.}
\textbf{465} (2009) 649--667], a new method for simulating nearly
linear stochastic partial differential equations (SPDEs) with additive
noise has been introduced. The key idea was to use suitable linear
functionals of the noise process in the numerical scheme which allow a
higher approximation order to be obtained. Following this approach, a
new simplified version of the scheme in the above named reference is
proposed and analyzed in this article. The main advantage of the
convergence result given here is the higher convergence order for
nonlinear parabolic SPDEs with additive noise, although the used
numerical scheme is very simple to simulate and implement.
\end{abstract}

%

\end{frontmatter}

\section{Introduction}
In this article, the numerical approximation of nonlinear parabolic
stochastic partial differential equations (SPDEs) is considered.
Following the idea in~\cite{jk09b} for somewhat linear SPDEs, a new
numerical method for simulating nonlinear SPDEs with additive noise is
proposed and analyzed in this article. The main advantage of the
convergence result in this article is the higher convergence order for
nonlinear parabolic SPDEs with additive noise in comparison to
convergence results of classical schemes such as the linear implicit
Euler scheme. Nevertheless, the here presented scheme is very simple
to simulate and implement.

More precisely, let
$T \in(0, \infty)$ be a real number, let
$( \Omega, \mathcal{F}, \mathbb{P})$ be a probability space
and let $H=L^2((0,1), \mathbb{R})$ be the $\mathbb{R}$-Hilbert\vadjust{\goodbreak}
space of equivalence classes of square integrable functions
from $(0,1)$ to $\mathbb{R}$.
Moreover, let $f\dvtx[0,1] \times\mathbb{R} \rightarrow
\mathbb{R}$ be a smooth function with bounded partial
derivatives, let $\xi\dvtx[0,1] \rightarrow\mathbb{R}$
with $\xi(0) = \xi(1) = 0$ be a smooth function and let
$W^Q \dvtx[0,T] \times\Omega\rightarrow H$ be
a standard $Q$-Wiener process
with a trace class operator
$Q\dvtx H \rightarrow H$ (see, e.g., Definition~2.1.9 in
\cite{pr07}).
It is a classical result (see, e.g., Proposition 2.1.5 in~\cite{pr07})
that the
covariance operator $Q\dvtx H \rightarrow H$ of the Wiener process
$W^Q \dvtx[0, T] \times\Omega\rightarrow H$
has an orthonormal basis $g_j \in H$, $j \in\mathbb{N}$,
of eigenfunctions with summable
eigenvalues $\mu_j \in[0,\infty)$, $j \in\mathbb{N}$.
In order to have a more concrete example,
we consider the choice
$g_j (x) = \sqrt{2} \sin(j \pi x)$
and $\mu_j =
c j^{-(r+1)}$ for all $x \in(0, 1)$ and
all $j \in\mathbb{N}$ with some
$ c \in[0,\infty)$ and some arbitrarily small
$ r \in(0,\infty)$
in the following and refer to Section~\ref{sec.assump} for our general setting.
Then we consider the SPDE
%
%
\begin{eqnarray} \label{eq:1}
dX_t &=& \biggl[ \frac{\partial^2}{\partial x^2} X_t + f(x, X_t) \biggr] \,dt +
dW_t^Q,\nonumber\\[-8pt]\\[-8pt]
X_t(0) &=& X_t(1) = 0,\qquad X_0 = \xi,\nonumber
\end{eqnarray}
for $x \in(0,1)$ and $t \in[0,T]$.
Under the assumptions above, the SPDE~(\ref{eq:1}) has a unique mild solution.
Specifically, there exists an up to indistinguishability unique
stochastic process $ X \dvtx[ 0, T ] \times\Omega\rightarrow H $ with
continuous sample paths which satisfies
%
%
\begin{equation}
X_t = e^{ A t } \xi+ \int^{ t }_{ 0 } e^{ A (t-s) } F( X_s ) \,ds +
\int^{ t }_{ 0 } e^{ A(t-s) } \,dW_s^Q ,\qquad\mathbb{P}\mbox{-a.s.}
\end{equation}
for all $t \in[0,T]$ where
$A\dvtx D(A) \subset H \rightarrow H$ is the Laplacian
with Dirichlet boundary conditions on $(0,1)$ and where
$F\dvtx H \rightarrow H $ is the
Nemytskii operator $( F(v) )(x) := f(x, v(x))$ for all
$x \in(0,1)$ and all $ v \in H$.

Then our goal is to solve the strong approximation problem of the
SPDE~(\ref{eq:1}). More precisely, we want to compute a
$ \mathcal{F} / \mathcal{B}(H) $-measurable
numerical approximation
$Y\dvtx\Omega\rightarrow H$ such that
%
%
\begin{equation} \label{eq:2}
\biggl( \mathbb{E} \biggl[ \int_0^1 | X_T (x) - Y(x) |^2 \,dx \biggr] \biggr)^{1/2} <
\varepsilon
\end{equation}
holds for a given precision $\varepsilon> 0$ with the least possible
computational effort (number of computational operations and
independent standard normal
random variables
needed to compute $Y \dvtx\Omega\rightarrow H$).
A computational operation
is here an arithmetical operation (addition, subtraction, multiplication,
division), a trigonometrical operation (sine, cosine) or an evaluation of
$f \dvtx(0, 1) \times\mathbb{R} \rightarrow\mathbb{R} $
or the exponential function.

In order to be able to calculate such a
numerical approximation on a computer, both the time interval
$[0,T]$ and the infinite-dimensional $\mathbb{R}$-Hilbert
space $H = L^2((0,1), \mathbb{R})$ have to be discretized.
While for temporal discretizations
the linear implicit Euler scheme is often used,
spatial discretizations are usually achieved with finite
elements, finite differences and spectral Galerkin methods.
For instance, the linear implicit Euler scheme combined with
spectral Galerkin methods which we denote by
$ \mathcal{F} /\mathcal{B}(H) $-measurable mappings
$Z_n^N \dvtx\Omega\rightarrow H$,
$n \in\{0, 1, \ldots, N^2 \}$, $N \in\mathbb{N}:=\{1,2,\ldots\}$,
is given by $Z_0^N := P_N( \xi)$
and
%
%
\begin{eqnarray} \label{eq:3}
Z_{n+1}^N &:=& \biggl( I - \frac{T}{N^2} A \biggr)^{ -1} \nonumber\\[-8pt]\\[-8pt]
&&{}\times\biggl( Z_n^N + \frac{T}{N^2}
\cdot(P_N F) ( Z_n^N ) + P_N \bigl( W^Q_{(n+1)T/N^2} -
W^Q_{nT/N^2} \bigr) \biggr)\nonumber
\end{eqnarray}
for every $n \in\{0, 1, \ldots, N^2 - 1 \}$ and every
$N \in\mathbb{N}$ where the bounded linear operators
$ P_N \dvtx H \rightarrow H $, $ N \in\mathbb{N} $, are given by
%
%
\begin{equation}
( P_N( v ) )( x ) := \sum_{ n = 1 }^N 2 \sin( n \pi x ) \int_0^1 \sin(
n \pi s ) v( s ) \,ds
\end{equation}
for all $ x \in( 0, 1 ) $, $ v \in H $ and all $ N \in\mathbb{N} $.
Note that the infinite-dimensional $ \mathbb{R} $-Hilbert space $ H $
is projected down to the $ N $-dimensional $ \mathbb{R} $-Hilbert space
$ P_N( H ) $ for the spatial discretization
and the time interval $ [0, T] $ is divided into $ N^2 $ subintervals,
that is, $ N^2 $ time steps are used, for the temporal discretization
in the scheme $ Z^N_n $,
$ n \in\{ 0,1,\ldots,N^2 \} $, above for $ N \in\mathbb
{ N
} $. The exact solution
$ X \dvtx[ 0, T ] \times\Omega\rightarrow H $
of the SPDE~(\ref{eq:1}) enjoys at least twice the
regularity in space than in time and
therefore, the quadratic number of time steps is used in
the scheme~(\ref{eq:3}) above (see also Walsh~\cite{w05b}
for details).

We now review how efficiently the numerical method~(\ref{eq:3})
solves the strong approximation problem~(\ref{eq:2}) of
the SPDE~(\ref{eq:1}).
Standard results in the literature (see, e.g.,
Theorem 2.1 in Hausenblas~\cite{h03a})
yield the existence of a real number $C> 0$
such that
%
%
\begin{equation} \label{eq:4}
\biggl( \mathbb{E} \biggl[ \int_0^1 | X_T (x) - Z_{N^2}^N(x) |^2 \,dx \biggr] \biggr)^{1/2}
\leq C\cdot N^{- 1}
\end{equation}
holds for all $N \in\mathbb{N}$.
Since $ P_N(H) $ is $N$-dimensional and since
$ N^2 $ time steps are used in~(\ref{eq:3}),
$O( N^3 \log(N))$
computational operations and independent standard normal
random variables
are needed
to compute $Z_{N^2}^N$ for $N \in\mathbb{N}$.
The log term in $ O( N^3 \log(N) ) $
for $ N \in\mathbb{N}$ arises due to
computing the nonlinearity with fast Fourier transform
(aliasing errors are neglected here).
Combining the computational effort $O( N^3 \log(N))$ and the estimate
(\ref{eq:4}) shows that the linear implicit Euler scheme needs
about $O( \varepsilon^{-3})$
computational operations and independent standard normal
random variables
to achieve a precision of size $\varepsilon> 0$ in the sense of~(\ref{eq:2}).
In fact, we have demonstrated that the linear implicit Euler scheme
method~(\ref{eq:3}) needs $O ( \varepsilon^{-(3 + \delta)} )$
computational operations and random variables
to solve~(\ref{eq:2}) for every arbitrarily small $\delta\in(0,
\infty)$
but for simplicity we write about $O ( \varepsilon^{-3} )$
computational operations and random variables here and below.

Recently, in~\cite{jk09b}, a new numerical method for simulating somewhat
linear SPDEs with additive noise has been introduced. The key
idea in~\cite{jk09b} is to use suitable linear functionals of the noise
process in the numerical scheme which allows a higher approximation
order to be obtained.
In this paper, we extend this idea to the case of nonlinear
SPDEs of the form~(\ref{eq:1}). More precisely, we introduce
the following numerical scheme which is a simplified version of the
scheme considered in~\cite{jk09b}. Let
$ Y^{N}_n \dvtx\Omega\rightarrow H $, $n \in\{0, 1, \ldots, N\}$, $N
\in
\mathbb{N}$, be
$ \mathcal{F} / \mathcal{B}(H) $-measurable
mappings given by
$ Y^{N}_0 := P_N( \xi) $
and
%
%
\begin{eqnarray} \label{eq:5}
Y^{N}_{n+1} &:=& e^{ A {T/N} }\biggl( Y^{N}_n + \frac{T}{N} \cdot(P_N
F)( Y^{N}_n ) \biggr) \nonumber\\[-8pt]\\[-8pt]
&&{} + P_N \biggl( \int_{nT/N}^{(n+1)T/N} e^{A
({(n+1)T}/{N} -s ) } \,dW_s^Q \biggr)\nonumber
\end{eqnarray}
$\mathbb{P}$-a.s. for every
$ n \in\{ 0, 1,\ldots, N-1 \} $ and every
$ N \in\mathbb{N} $.
Note that the infinite-dimensional $ \mathbb{R}$-Hilbert space $ H $ is
projected down to the $ N $-dimensional $ \mathbb{R} $-Hilbert space $
P_N( H ) $ for the spatial discretization
and the time interval $ [0, T] $ is divided into $ N $ subintervals,
that is, $ N $ time steps are used, for the temporal discretization
in the scheme $ Y^N_n $, $ n \in
\{ 0,1,\ldots,N \} $, above for $ N \in\mathbb{ N } $.

We now illustrate the main result of this article
(Theorem~\ref{theorem1})
and show how efficiently the
method~(\ref{eq:5}) solves
the strong approximation problem~(\ref{eq:2})
of the SPDE~(\ref{eq:1}).
Theorem~\ref{theorem1} shows the existence of real numbers
$C_\delta> 0$, $\delta\in(0,1)$, such that
%
%
\begin{equation} \label{eq:6}
\biggl( \mathbb{E} \biggl[ \int_0^1 | X_T (x) - Y_{N}^N(x) |^2 \,dx \biggr] \biggr)^{1/2}
\leq C_\delta\cdot N^{(\delta- 1)}
\end{equation}
holds for all $N \in\mathbb{N}$ and all arbitrarily small $\delta\in(0,1)$.
The stochastic integrals
%
%
\begin{equation}\label{sint}
P_N \biggl( \int_{nT/N}^{(n+1)T/N} e^{A ({(n+1)T}/{N} -s
) } \,dW_s^Q \biggr)
\end{equation}
for $n \in\{0, 1, \ldots, N\}$
and $N \in\mathbb{N}$ in~(\ref{eq:5})
provide more information
about the exact solution and this allows us
to obtain the estimate~(\ref{eq:6}) although
only $N$ time steps (instead of $N^2$ time steps in the case
of the linear implicit Euler scheme)
are used in~(\ref{eq:5}).
Nevertheless, since the stochastic integrals~(\ref{sint})
in~(\ref{eq:5}) depend linearly on the Wiener
process $W^Q \dvtx[0,T] \times\Omega\rightarrow H$,
they are again normally distributed and hence easy to
simulate.
More precisely, since $ P_N(H) $ is $N$-dimensional and since
$ N $ time steps are used in~(\ref{eq:5}),
$O( N^2 \log(N))$
computational operations and independent standard normal
random variables\vadjust{\goodbreak}
are needed
to compute $Y_{N}^N$ for $N \in\mathbb{N}$.
The log term in $ O( N^2 \log(N) ) $
for $ N \in\mathbb{N}$ also arises due to
computing the nonlinearity with fast Fourier
transform (aliasing errors are neglected here).
Combining the computational effort $O( N^2 \log(N))$ and the estimate
(\ref{eq:6}) shows that the numerical scheme~(\ref{eq:5}) needs
about $O( \varepsilon^{-2})$
computational operations and independent standard normal
random variables
to achieve a precision of size $\varepsilon> 0$ in the sense of~(\ref{eq:2}).

The estimates~(\ref{eq:4}) and~(\ref{eq:6})
are both asymptotic results
since there is no information
about the size of the corresponding error
constants. 
In particular, the error constants $C_\delta\in(0, \infty) $,
$\delta
\in(0,1)$,
in~(\ref{eq:6}) could be much bigger than in~(\ref{eq:4}).
Therefore, from a practical point of view,
one may ask whether the numerical method~(\ref{eq:5})
solves the strong approximation problem~(\ref{eq:2})
more efficiently
than the linear implicit Euler scheme~(\ref{eq:3})
for a given example of the form~(\ref{eq:1})
and a given concrete $ \varepsilon> 0 $.
In order to analyze this question, we
compare both methods in the case of a simple
reaction diffusion SPDE of the
form~(\ref{eq:1})
(see Section~\ref{example1}
for details) and assume that the strong
approximation problem~(\ref{eq:2}) should
be solved with the precision $ \varepsilon= \frac{1}{300} $.
In that example, it turns out
that the linear implicit Euler scheme precisely needs
$2^{21} = 2\mbox{,} 097\mbox{,}152$
independent standard normal random variables while
the numerical method~(\ref{eq:5})
precisely needs $2^{16} = 65\mbox{,} 536$ independent standard
normal random variables to achieve an approximation error of size
$\varepsilon= \frac{1}{300}$
(see Tables~\ref{tab1} and~\ref{tab2} in Section~\ref{example1}).
We also emphasize that the numerical
scheme~(\ref{eq:5}) is very simple to implement and refer
to Figure~\ref{code} for a short \textsc{matlab} code.

Having illustrated the main result of this article,
we now sketch the key idea in the proof
of Theorem~\ref{theorem1}.
The main difficulty was to estimate the discretization error
for nonlinear $F$.
In that case, the main problem was to establish
estimates of the form
%
%
\begin{equation}\label{eq:toshow2}
\Biggl\| \sum_{n = 0}^{N - 1} \int_{ { n T }/{ N } }^{ { ( n + 1 ) T
}/{ N } } e^{ A { ( T - s ) } } \bigl( F(X_s) - F( X_{ { n T }/{ N } } )
\bigr) \,ds \Biggr\|_{ L^2 ( \Omega; H ) } \leq C_\delta\cdot
N^{(\delta-1)}\hspace*{-32pt}
\end{equation}
for all $ N \in\mathbb{N} $ and all $\delta\in(0,1)$ where
$C_\delta\in(0, \infty) $, $\delta\in(0,1)$, are appropriate constants
and where we write $ \| Y \|_{ L^2 ( \Omega; H ) } := ( \mathbb{E} [
\int_0^1 | Y(x) |^2 \,dx ] )^{ 1/2 } \in[0,\infty] $ for every $
\mathcal{F} / \mathcal{B}(H) $-measurable mapping $ Y
\dvtx\Omega\rightarrow H $ for simplicity. The smoothness of the
Nemytskii operator $ F $ on an appropriate subspace $V \subset H$
shows that it remains to estimate
%
%
\begin{eqnarray}\label{eq:toshow}
&&\Biggl\| \sum_{n = 0}^{N - 1} \int_{ n T / N }^{ { ( n + 1 ) T
}/{ N } } e^{ A { ( T - s ) } } F'( X_{ n T / N } ) ( X_s -
X_{ n T / N } ) \,ds \Biggr\|_{ L^2 ( \Omega; H ) } \nonumber\\[-8pt]\\[-8pt]
&&\qquad\leq
C_\delta\cdot N^{(\delta- 1)}\nonumber
\end{eqnarray}
for all $ N \in\mathbb{N}$ and all $\delta\in(0,1)$.
In~\cite{jk09b}, the linear operators
$ F'(v) $ for $ v \in H $
and $ A \dvtx D(A) \subset H
\rightarrow H $ are assumed
to commute in some sense
which is fulfilled in the case of linear\vadjust{\goodbreak} $F$
such as $ F(v) = v $, $ v \in H $,
but excludes nonlinear
Nemytskii operators such as
$ F(v) = \frac{ ( 1 - v ) }{ ( 1 + v^2 ) } $,
$ v \in H $ (see Assumption 2.4
in~\cite{jk09b} for details).
Under this commutativity condition,~(\ref{eq:toshow})
can easily be established by using the smoothing
effect of the semigroup $ e^{ A t } $,
$ t \in[0,T] $ (see Section 5.b.i
in~\cite{jk09b}).
Instead of this condition, our key assumption
on the nonlinearity is an appropriate
estimate on the adjoint operators
of the Fr\'{e}chet derivative operators
of $F$
[see~(\ref{assumpt2_2})].
Since in our examples $F$ is
a (nonlinear) Nemytskii operator, the
derivative operators $ F'(v) $, $v \in V $, are
self-adjoint and hence, it can easily be seen that this assumption
is fulfilled [see~(\ref{eqRef1}) in
Section~\ref{secex} for details].
Moreover, this assumption enables use to
show~(\ref{eq:toshow}) and hence~(\ref{eq:toshow2})
[see Section~\ref{temp_error}
and particularly estimate~(\ref{keyproof})].
We also mention that the difficulty to estimate
(\ref{eq:toshow2}) can be avoided
by using a more complicated scheme with a second
linear functional (see Section 6.4 in~\cite{jk09c}).

Finally, we would like to point out limitations of the here
presented numerical method. The following assumption
is essential to apply our algorithm. The eigenfunctions
of the dominating linear operator and of the covariance
operator of the driving additive noise process of
the SPDE must coincide and must be known explicitly.

The rest of this article is organized as follows.
The basic setting and the assumptions that we use
(including our key assumption on the adjoint
of the Fr\'{e}chet derivative of the nonlinearity)
are presented in Section~\ref{sec.assump}.
The new numerical scheme and its convergence
theorem which is the main result of
this article are given in Section~\ref{main}.
This result is illustrated with some examples
and some numerical simulations in Section~\ref{secex}.
Although our setting in Section~\ref{sec.assump} uses
the standard global Lipschitz assumption on the nonlinearity
of the SPDE, we demonstrate the efficiency of our method numerically
for a SPDE with a cubic nonglobally Lipschitz nonlinearity in
Section~\ref{sec.furtherexamples}.
Proofs are postponed to the final section.
%

\section{Setting and assumptions}\label{sec.assump}

Fix $ T \in(0, \infty) $ and let
$ ( \Omega, \mathcal{F}, \mathbb{P} ) $
be a probability space with a normal filtration
$( \mathcal{F}_t )_{t \in[ 0, T] }$ which means
$\mathcal{F}_{t+} = \mathcal{F}_t$ for all $t \in[0,T)$
and
$\{ A \in\mathcal{F} | \mathbb{P}[A] = 0 \} \subset\mathcal{F}_0$
(see, e.g., Definition 2.1.11 in~\cite{pr07}).
In addition, let
$ ( V, \| \cdot\|_V ) $
be a separable $ \mathbb{R} $-Banach space
and let
$( H, \langle \cdot, \cdot\rangle_H
, \| \cdot\|_H)$ be a separable
$\mathbb{R}$-Hilbert space
with $ V \subset H $ continuously.
The following assumptions will be used.
\begin{assumption}[(Linear operator $A$)] \label{semigroup}
Let $( \lambda_n)_{n \in\mathbb{N} } \subset(0, \infty) $
be an increasing sequence of real numbers and let
$ (e_n)_{ n \in\mathbb{N} } \subset H$
be an orthonormal basis of $H$.
Assume that the linear operator
$ A \dvtx D(A) \subset H \rightarrow H $
is given by
\[
Av = \sum_{n=1}^\infty- \lambda_n \langle e_n, v
\rangle_H
e_n
\]
for all $v \in D(A)$ with $D(A) = \{ w \in H | {\sum_{n=1}^\infty}|
\lambda_n |^2 | \langle e_n, w \rangle_H |^2 < \infty\}$.\vadjust{\goodbreak}
\end{assumption}

Let $ D( (-A)^r ) $ with $ \| v \|_{ D( (-A)^r ) } = \| (-A)^r v \|_H $
for $ v \in D( (-A)^r ) $ and $ r \in\mathbb{R} $ denote the domains of
fractional powers of the linear operator $ - A $ (see, e.g., Section
3.7 in~\cite{sy02}).
\begin{assumption}[(Nonlinearity $F$)] \label{drift}
Assume $ D( (-A)^{1/2} ) \subset V $ continuously and let $ F
\dvtx V
\rightarrow V $ be a twice continuously Fr\'{e}chet differentiable
mapping with
%
%
\begin{eqnarray}
\label{assumpt2_1}
\| F'(v) w \|_H &\leq& c \| w \|_H,\nonumber\\[-8pt]\\[-8pt]
\| F'(v)\|_{ L ( V ) } &\leq& c ,\qquad \|F''(v) \|_{ L^{(2)} (V) } \leq c
,\nonumber
\\
\label{assumpt2_2}
\|
( F'(u) )^* 
\|_{L ( D ( (-A)^{1/2} ) ) } &\leq& c \bigl( 1 + \| u \|_{ D (
(-A)^{1/2} ) } \bigr)
\end{eqnarray}
for every $ v, w \in V $ and every
$u \in{ D ( (-A)^{1/2} ) }$
where $ c \in[0,\infty) $
is a given real number.
\end{assumption}

By definition $ F'(v) \in L(V) $ is a bounded linear mapping from $ V $
to $ V $ for every $ v \in V $. Due to the first condition in (\ref
{assumpt2_1}), we also have that
$ F'(v) \in L(H) $ is a bounded linear mapping from $ H $ to $ H $ for every $ v \in V $. In
that sense,
the adjoint operator $ ( F'(v) )^* \in L(H) $ given by
\[
\langle( F'(v) )^* u, w \rangle_H = \langle u, F'(v) w \rangle_H
\]
for all $ u, w \in H $ is well defined for every $ v \in V $. Due to
(\ref{assumpt2_2}), the operator
$( F'(v) )^* \in L ( H )$
is also a bounded linear mapping from
${ D ( (-A)^{1/2} ) }$ to
${ D ( (-A)^{1/2} ) }$ for every
$v \in{ D ( (-A)^{1/2} ) }$.
\begin{assumption}[(Stochastic process $O$)] \label{stochconv}
Let
$ O \dvtx[0,T] \times\Omega\rightarrow
D ( (-A)^\gamma) $
be a centered and adapted stochastic process with
continuous sample paths such that
$O_ {t_2} - e^{ A {( t_2 - t_1 )} } O_{t_1}$
is independent of $\mathcal{F}_{t_1}$ for all
$0 \leq t_1 < t_2 \leq T$ and such that
\[
\mathbb{E} \Bigl[ \sup_{0 \leq t \leq T} \| (-A)^\gamma O_t \|_H^4 \Bigr] +
\sup_{0 \leq t_1 < t_2 \leq T} \bigl( ( t_2 - t_1 )^{ - 4 \theta}
\mathbb{E}[ \| O_{t_2} - O_{t_1} \|_V^4 ] \bigr) < \infty
\]
holds where $ \gamma\in[ \frac{1}{2}, 1 ) $ and $ \theta\in(0,
\frac{1}{2} ]$ are given real numbers.
\end{assumption}
\begin{assumption}[(Initial value $\xi$)] \label{initial}
Let $ \xi\dvtx\Omega\rightarrow D(A) $
be a $ \mathcal{F}_0 /\mathcal{B}(D(A)) $-measurable mapping
with
$
\mathbb{E}
[
\|
A \xi
\|_H^4
]
<
\infty
$.
\end{assumption}

%
%
%
%
%
These assumptions suffice to ensure the existence of a
unique solution of the SPDE~(\ref{eqSolution}).
\begin{lemma}[(Existence of the solution)] \label{existence}
Let Assumptions~\ref{semigroup}--\ref{initial} be fulfilled.
Then there exists a unique adapted stochastic process
$X \dvtx[0,T] \times\Omega\rightarrow D( (-A)^\gamma)$
with continuous sample paths which fulfills
%
%
\begin{equation}\label{eqSolution}
X_t(\omega) = e^{At} \xi(\omega) + \int_0^t e^{A(t-s)} F(X_s(\omega))
\,ds + O_t(\omega)\vadjust{\goodbreak}
\end{equation}
for all $t \in[0,T]$ and all $\omega\in\Omega$. Moreover, $X\dvtx[0,T]
\times\Omega\rightarrow D( (-A)^\gamma)$ satisfies $ \mathbb{E}
[{\sup_{0 \leq t \leq T}} \|(-A)^\gamma X_t\|_H^4] < \infty $.
\end{lemma}

The proof of Lemma~\ref{existence} is given in Section~\ref{proofs}.
Some examples satisfying Assumptions~\ref{semigroup}--\ref{initial} are
presented in Section~\ref{secex}.

\section{Numerical scheme and main result} \label{main}

For numerical approximations of the SPDE~(\ref{eqSolution}), we have to
discretize both the time interval $[0,T]$ and the $\mathbb{R}$-Hilbert
space $H$. To this end, we use projections $P_N\dvtx H \rightarrow H$
given by $P_N (v) := \sum_{n = 1}^N \langle e_n, v \rangle_H e_n$ for
every $v \in H$, $N \in\mathbb{N}$ and finite-dimensional
$\mathbb{R}$-Hilbert spaces $ H_N \subset H $ given by $H_N := P_N (H)$
for every $N \in\mathbb{N}$. Finally, we define $ \mathcal{F} /
\mathcal{B}(H_N) $-measurable mappings $ Y^{N,M}_m
\dvtx\Omega\rightarrow H_N $ for $ m \in\{ 0, 1,\ldots, M \} $ and $ N,
M \in\mathbb{N} $ by $ Y^{N,M}_0(\omega) := P_N( \xi(\omega) ) + P_N(
O_0(\omega) ) $ and by
%
%
\begin{eqnarray} \label{expEulerScheme}
Y^{N,M}_{m+1}(\omega) &:=& e^{ A {T/M} }\biggl( Y^{N,M}_m(\omega) +
\frac{T}{M} \cdot(P_N F)( Y^{N,M}_m(\omega) ) \biggr)
\nonumber\\[-8pt]\\[-8pt] 
&&{} + P_N\bigl( O_{ { (m+1) T }/{ M } }(\omega) - e^{ A {T/M} }
O_{ { m T }/{ M } }(\omega) \bigr)\nonumber
\end{eqnarray}
for every
$ m \in\{ 0, 1,\ldots, M-1 \} $,
$ N, M \in\mathbb{N} $ and every
$ \omega\in\Omega$.
In many examples, this scheme is as easy to
simulate as the classical linear
implicit Euler scheme.
We refer to Section~\ref{secex}
for a detailed description of the implementation of
our numerical scheme including a short \textsc{matlab} code.
\begin{theorem}\label{theorem1}
Let Assumptions~\ref{semigroup}--\ref{initial}
be fulfilled. Then there is
a real number $C > 0$ such that
%
%
\begin{equation} \label{main_equation}
\bigl( \mathbb{E} \bigl[ \| X_{ {mT}/{M} } - Y^{N,M}_m \|_H^2 \bigr] \bigr)^{1/2}
\leq C \biggl( \frac{ 1 }{ (\lambda_N)^\gamma} + \frac{ ( 1 + \log( M ) ) }{
M^{ 2 \theta} } \biggr)
\end{equation}
holds for every $m \in\{ 0, 1,\ldots, M \}$ and every $ N, M \in
\mathbb{N} $ where $ (\lambda_N)_{N \in\mathbb{N}} \subset(0, \infty)$
is given in Assumption~\ref{semigroup} where $ \gamma\in[ \frac{1}{2},
1 )$ and $ \theta\in(0, \frac{1}{2}] $ are given in
Assumption~\ref{stochconv} where $ X \dvtx[0,T] \times\Omega\rightarrow
D( ( - A )^\gamma) $ is the solution of the SPDE~(\ref{eqSolution}) and
where $ Y^{N,M}_m \dvtx\Omega\rightarrow H_N $, $ m \in\{ 0,1,\ldots, M
\} $, $ N, M \in\mathbb{N} $, is given by~(\ref{expEulerScheme}).
\end{theorem}

Here and below $ \log$ is the natural logarithm. While the expression
$ \frac{ 1 }{ (\lambda_N)^\gamma}$ for $N \in\mathbb{N} $ in
(\ref{main_equation}) arises due to discretizing the infinite-dimensional
$\mathbb{R}$-Hilbert space~$H$, the expression
$\frac{ ( 1 + \log( M ) ) }{M^{ 2 \theta}} $
for $M \in\mathbb{N}$ arises due to discretizing the time interval $[0,T]$.
We would like to remark that the logarithmic term in
$\frac{ ( 1 + \log( M ) ) }{M^{ 2 \theta}} $ for $M \in\mathbb{N}$ can
be avoided by assuming
$F(D((-A)^{1/2})) \subset D((-A)^\varepsilon) $ and an
appropriate linear growth condition on $F$ for some $\varepsilon> 0$.
Although this condition is fulfilled in our examples below, we use this
logarithmic term in Theorem~\ref{theorem1} here in order to
formulate Assumption~\ref{drift} in our abstract setting as
simple as possible.

A similar result could be obtained for SPDEs of the
form~(\ref{eqSolution}) but with a time dependent nonlinearity $F$.
However, we omit the time dependency of the nonlinearity here for
simplicity.

\section{Examples} \label{secex}

Let $H = L^2((0,1), \mathbb{R} )$ be the
$\mathbb{R}$-Hilbert space of equivalence classes of
$\mathcal{B} ((0,1))/\mathcal{B}(\mathbb{R})$-measurable
and square integrable functions from $(0,1)$ to $\mathbb{R}$
with the scalar product and the norm given by
\[
\langle v, w \rangle_H = \int_0^1 v(s) w(s) \,ds,\qquad \| v \|_H = \biggl( \int_0^1
| v(s) |^2 \,ds \biggr)^{1/2}
\]
for every $v, w \in H$.
In addition, let $V = C( [0,1], \mathbb{R}) $ be the
$\mathbb{R}$-Banach space of continuous functions from
$[0,1]$ to $\mathbb{R}$ equipped with the norm
$
\| v \|_V = {\sup_{ 0 \leq x \leq1}}
| v(x) |
$
for every $v \in V$.

Let $\kappa\in(0, \infty)$ be a given
positive real number and let
$(\lambda_n)_{n \in\mathbb{N}} \subset(0, \infty)$
and $( e_n )_{n \in\mathbb{N} } \subset H$
be given by
\[
\lambda_n := \kappa n^2 \pi^2,\qquad e_n(x) := \sqrt{2} \sin(n \pi x)
\]
for every $x \in(0,1)$ and every $n \in\mathbb{N}$.
Hence, the linear operator
$A\dvtx D(A) \subset H \rightarrow H$
reduces to the Laplacian with Dirichlet boundary conditions
on the interval $(0,1)$ times
the constant $ \kappa \in (0,\infty)$ (see, e.g., Section 3.8.1
in Sell and You~\cite{sy02}).
In particular, ${ D ( (-A)^{1/2} ) }$
reduces to the $\mathbb{R}$-Sobolev space
$H_0^1( (0,1), \mathbb{R} )$ equipped with the norm
\begin{eqnarray*}
\| u \|_{ D ( (-A)^{1/2} ) } & = & \| (-A)^{1/2} u \|_H \\
&=& \Biggl(
\sum_{n=1}^\infty\kappa n^2 \pi^2 | \langle e_n, u \rangle_H |^2
\Biggr)^{1/2}
\\
&=& \sqrt{ \kappa} \biggl( \int_0^1 | u'(x) |^2 \,dx \biggr)^{1/2}
\end{eqnarray*}
for all $u \in{ D ( (-A)^{1/2} ) }$.
(See Sell and You~\cite{sy02} for more information about this space.)

Furthermore, let $f\dvtx[0,1] \times\mathbb{R} \rightarrow\mathbb{R}$
be a twice continuously differentiable function with the bounded
partial derivatives
\[
\biggl| \biggl( \frac{\partial f}{\partial y} \biggr) (x, y) \biggr| \leq K,\qquad
\biggl| \biggl(\frac{\partial^2 f}{\partial x\, \partial y} \biggr) (x, y) \biggr| \leq K,\qquad
\biggl| \biggl(\frac{\partial^2 f}{\partial y^2} \biggr) (x, y) \biggr| \leq K
\]
for all $ x \in[0,1]$ and all $ y \in\mathbb{R}$
with an arbitrary constant $K \in[0,\infty)$. Then
the Nemytskii operator $F\dvtx V \rightarrow V$ given by
$(F(v))(x) = f(x, v(x))$ for every $ x \in[0,1]$
and every $v \in V$ satisfies Assumption~\ref{drift}.
To see this note that
\begin{eqnarray*}
F'(u)(v) &=& \biggl( \frac{ \partial f}{ \partial y} \biggr) (x,u(x)) \cdot v(x),
\\
F''(u)(v,w) &=& \biggl( \frac{ \partial^2 f}{ \partial y^2 } \biggr) (x,u(x)) \cdot
v(x) \cdot w(x)
\end{eqnarray*}
holds for all $u, v, w \in V$. Therefore, we have
\begin{eqnarray*}
&& \| ( F'(u) )^* v \|_{ D ( (-A)^{1/2} ) }^2
\\[-1pt]
&&\qquad
= \| (-A)^{1/2} F'(u) v \|_H^2\\[-1pt]
&&\qquad= \kappa\int_0^1 \biggl|
\frac{\partial}{\partial x} \biggl\{ \biggl( \frac{\partial f}{\partial y} \biggr) (x,
u(x)) \cdot v(x) \biggr\} \biggr|^2 \,dx
\\[-1pt]
&&\qquad
= \kappa\int_0^1 \biggl| \biggl\{ \frac{\partial}{\partial x} \biggl[ \biggl( \frac{\partial
f}{\partial y} \biggr) (x, u(x)) \biggr] \biggr\} v(x) + \biggl( \frac{\partial f}{\partial y}
\biggr) (x, u(x)) \cdot v'(x) \biggr|^2 \,dx
\\[-1pt]
&&\qquad\leq
2 \kappa\int_0^1 \biggl| \biggl\{ \frac{\partial}{\partial x} \biggl[ \biggl( \frac{\partial
f}{\partial y} \biggr) (x, u(x)) \biggr] \biggr\} v(x) \biggr|^2 \,dx\\[-1pt]
&&\qquad\quad{} + 2 \kappa\int_0^1 \biggl| \biggl(
\frac{\partial f}{\partial y} \biggr) (x, u(x)) \cdot v'(x) \biggr|^2 \,dx
\end{eqnarray*}
and
\begin{eqnarray*}
&& \| ( F'(u) )^* v \|_{ D ( (-A)^{1/2} ) }^2
\\[-1pt]
&&\qquad\leq
2 \kappa\| v \|_V^2 \int_0^1 \biggl| \frac{\partial}{\partial x} \biggl[ \biggl(
\frac{\partial f}{\partial y} \biggr) (x, u(x)) \biggr] \biggr|^2 \,dx\\[-1pt]
&&\qquad\quad{} + 2 \kappa K^2
\int_0^1 | v'(x) |^2 \,dx
\\[-1pt]
&&\qquad\leq
4 \| v \|_{ D ( (-A)^{1/2} ) }^2 \int_0^1 \biggl| \biggl( \frac{\partial^2
f}{\partial x \,\partial y} \biggr) (x, u(x)) \biggr|^2 \,dx
\\[-1pt]
&&\qquad\quad{}
+ 4 \| v \|_{ D ( (-A)^{1/2} ) }^2 \int_0^1 \biggl| \biggl( \frac{\partial^2
f}{\partial y^2} \biggr) (x, u(x))\cdot u'(x) \biggr|^2 \,dx\\[-1pt]
&&\qquad\quad{} + 2 K^2 \| v \|_{ D (
(-A)^{1/2} ) }^2
\end{eqnarray*}
for all $u, v \in{ D ( (-A)^{1/2} ) }$.
Hence, we obtain
\begin{eqnarray*}
&& \| ( F'(u) )^* v \|_{ D ( (-A)^{1/2} ) }^2
\\[-1pt]
&&\qquad\leq
4 K^2 \| v \|_{ D ( (-A)^{1/2} ) }^2\\[-1pt]
&&\qquad\quad{} + 4 K^2 \| v \|_{ D (
(-A)^{1/2} ) }^2 \int_0^1 | u'(x) |^2 \,dx\\[-1pt]
&&\qquad\quad{} + 2 K^2 \| v \|_{ D (
(-A)^{1/2} ) }^2
\\[-1pt]
&&\qquad
= 6 K^2 \| v \|_{ D ( (-A)^{1/2} ) }^2 + 4 K^2 \kappa^{-1} \| v
\|_{ D ( (-A)^{1/2} ) }^2 \| u \|_{ D ( (-A)^{1/2} ) }^2
\end{eqnarray*}
and
%
%
\begin{eqnarray} \label{eqRef1}
&& \| ( F'(u) )^* v \|_{ D ( (-A)^{1/2} ) }
\nonumber\\
&&\qquad\leq
\sqrt{ 6 K^2 \| v \|_{ D ( (-A)^{1/2} ) }^2 + 4 K^2 \kappa^{-1}
\| v \|_{ D ( (-A)^{1/2} ) }^2 \| u \|_{ D ( (-A)^{1/2} )
}^2 }\nonumber\\
&&\qquad\leq
\sqrt{6} K \| v \|_{ D ( (-A)^{1/2} ) } + 2 K \kappa^{ -1/2 } \| v \|_{ D ( (-A)^{1/2} ) } \| u \|_{ D(
(-A)^{1/2} ) }
\\
&&\qquad\leq
K \| v \|_{ D( (-A)^{1/2} ) } \bigl( 3 + 2 \kappa^{ -1/2 } \|
u \|_{ D( (-A)^{1/2} ) } \bigr) \nonumber\\
&&\qquad\leq
( 3 + 2 \kappa^{ -1/2 } ) K \| v \|_{ D (
(-A)^{1/2} ) } \bigl( 1 + \| u \|_{ D ( (-A)^{1/2} ) }
\bigr)\nonumber
\end{eqnarray}
for all $u, v \in{ D ( (-A)^{1/2} ) }$.
This shows that $F$ indeed satisfies Assumption~\ref{drift} with $c = 3K$.

Let $(b_n)_{n \in\mathbb{N}} \subset\mathbb{R}$ be a sequence of real
numbers with
$\sum_{n=1}^\infty n^\varepsilon|b_n|^2 < \infty$ for some
arbitrarily small $\varepsilon\in(0,\infty)$.
Lemma 4.3 in~\cite{bj09} then gives the existence
of an up to indistinguishability unique stochastic process
$O\dvtx[0,T] \times\Omega\rightarrow V$
which satisfies Assumption~\ref{stochconv} for
$\theta= \frac{1}{2}$ and $\gamma= \frac{1}{2}$ and which satisfies
\[
\mathbb{P} \Biggl[ O_t = \sum_{n=1}^\infty b_n \biggl( \int_0^t e^{ - \lambda_n (t
- s) } \,d \beta_s^n \biggr) e_n \Biggr] = 1
\]
for all $t \in[0,T]$ where the
$ \beta^n \dvtx[0,T] \times\Omega\rightarrow\mathbb{R} $,
$ n \in\mathbb{N}$, are independent standard Brownian
motions with respect to a given normal filtration
$( \mathcal{F}_t )_{t \in[0,T]}$.

Moreover, the mapping $ \xi\dvtx\Omega\rightarrow V $ given by
\[
(\xi(\omega) ) (x) = \tfrac{1}{ \sqrt{2} } \sin{ (\pi x) } + \tfrac{ 3
\sqrt{2} } {5} \sin{ ( 3 \pi x ) }
\]
for all $\omega\in\Omega$ and all $ x \in(0,1)$ obviously
satisfies Assumption~\ref{initial}.

In view of the above choice,
the SPDE~(\ref{eqSolution}) reduces to
%
%
\begin{eqnarray} \label{reducedForm}
dX_t &=& \biggl[ \kappa\,\frac{\partial^2}{\partial x^2} X_t + f(x, X_t) \biggr] \,dt +
B \,dW_t,
\nonumber\\
X_t(0) &=& X_t(1) = 0,\\
X_0(x) &=& \frac{ \sin(\pi x) }{ \sqrt{2}} + \frac{3
\sqrt{2}}{5} \sin( 3\pi x )\nonumber
\end{eqnarray}
for $x \in[0,1]$ and
$ t \in[0,T]$
where the linear operator $B\dvtx H \rightarrow H$
is given by
\[
Bv = \sum_{n=1}^\infty b_n \langle e_n, v \rangle e_n
\]
for all $v \in H$ and where
$( W_t )_{ t \in[0,T] }$ is a
cylindrical $I$-Wiener process on $H$.

Since Assumption~\ref{stochconv} is fulfilled
for
$\theta= \frac{1}{2}$ and $\gamma= \frac{1}{2}$,
Theorem~\ref{theorem1} shows the existence
of a real number $C > 0$,
such that
%
%
\begin{equation} \label{errorestimate2}\quad
\biggl( \mathbb{E} \biggl[ \int_0^1 | X_{nT/M} (x) - Y_n^{N,M} (x) |^2 \,dx \biggr]
\biggr)^{1/2} \leq C \biggl( \frac{ 1 }{ N } + \frac{ ( 1 + \log(M) )}{M} \biggr)\vadjust{\goodbreak}
\end{equation}
holds for all $n \in\{0, 1, \ldots, M\}$ and all $N, M \in\mathbb{N}$.
While the expression $\frac{1}{N}$ for $N \in\mathbb{N}$ in
(\ref{errorestimate2}) corresponds to the spatial discretization
error, the expression $\frac{ ( 1 + \log(M) )}{M} $ for
$M \in\mathbb{N}$ in~(\ref{errorestimate2}) corresponds to the
temporal discretization error. Since these error terms are
nearly of the same size, we choose $M=N$ and
consider the numerical approximations $Y_n^{N,N} \dvtx\Omega
\rightarrow H_N$,
$n \in\{0, 1, \ldots, N\}$, $N \in\mathbb{N}$, in the following.
Due to~(\ref{errorestimate2}), we obtain the existence
of real numbers $C_\delta> 0$, $\delta\in(0, 1 )$,
such that
%
%
\begin{equation} \label{particularForm}
\biggl( \mathbb{E} \biggl[ \int_0^1 | X_T (x) - Y_N^{N,N} (x) |^2 \,dx \biggr]
\biggr)^{1/2} \leq C_\delta\cdot N^{(\delta- 1)}
\end{equation}
holds for all $N \in\mathbb{N}$ and all arbitrarily small $\delta\in(0,1)$.

In order to describe the implementation of the
numerical scheme~(\ref{expEulerScheme}) in this example,
we use the
$ \mathcal{F} / \mathcal{B}(\mathbb{R}) $-measurable
mappings
$ Y^{N,M}_{n,m} \dvtx\Omega\rightarrow\mathbb{R}$
given by
$ Y^{N,M}_{n,m} (\omega) :=
\langle e_n, Y^{N,M}_m (\omega) \rangle_H$
for all $n \in\{ 1, 2, \ldots, N \}$,
$m \in\{ 0, 1, \ldots, M \}$ and all
$N, M \in\mathbb{N}$.
The numerical scheme
(\ref{expEulerScheme}) for the SPDE~(\ref{reducedForm}) with $M$ $=$ $N$
then reduces to
$ Y^{N,N}_{1,0} = \frac{1}{2} $,
$ Y^{N,N}_{2,0} = 0 $,
$ Y^{N,N}_{3,0} = \frac{3}{5} $,
$ Y^{N,N}_{4,0} = Y^{N,N}_{5,0} = \cdots= 0 $ and
%
%
\begin{eqnarray} \label{reducedForm2}
Y^{N,N}_{1,n+1} &=& e^{- \kappa\pi^2 {T}/{N}} \biggl( Y^{N,N}_{1,n} +
\frac{T}{N} \langle e_1, F ( Y^{N,N}_n ) \rangle_H \biggr)\nonumber\\
&&{} + \frac{\sqrt{1 -
e^{ - 2 \kappa\pi^2 {T}/{N} } } } {b_1 \cdot\pi\sqrt{2 \kappa}}
\chi_{1,n}^N, \nonumber\\
Y^{N,N}_{2,n+1} &=& e^{- \kappa\pi^2 2^2 {T/N}} \biggl( Y^{N,N}_{2,n} +
\frac{T}{N} \langle e_2, F ( Y^{N,N}_n ) \rangle_H \biggr)\nonumber\\
&&{} + \frac{ \sqrt{1
- e^{ - 2 \kappa2^2 \pi^2 {T/N} } } } {b_2 \cdot2 \pi\sqrt{2
\kappa}} \chi_{2,n}^N,
\\
\vdots
&&
\vdots\nonumber
\\
Y^{N,N}_{N,n+1} &=& e^{- \kappa\pi^2 N^2 {T/N}} \biggl( Y^{N,N}_{N,n} +
\frac{T}{N} \langle e_N, F ( Y^{N,N}_n ) \rangle_H \biggr)\nonumber\\
&&{} + \frac{ \sqrt{1
- e^{ - 2 \kappa N^2 \pi^2 {T/N} } } } {b_N \cdot N \pi\sqrt{2
\kappa}} \chi_{N,n}^N\nonumber
\end{eqnarray}
for all $n \in\{0, 1, \ldots, N-1 \}$ and all $N \in\mathbb{N}$
where the
$ \mathcal{F} / \mathcal{B}(\mathbb{R}) $-measurable
mappings
$ \chi^N_{n,m} \dvtx\Omega\rightarrow\mathbb{R}$ for
$n \in\{1, 2, \ldots, N\}$, $m \in\{0, 1, \ldots, N-1\}$ and
$N \in\mathbb{N}$ are independent standard normal random
variables. Since $O( N^2 \log(N))$
computational operations and independent standard normal
random variables
(computational
effort) are needed to compute the numerical solution
$Y_N^{N,N}$ given by~(\ref{reducedForm2}) for $N \in\mathbb{N}$,
it follows that $Y_N^{N,N}$ converges with order $\frac{1}{2}^-$ with
respect to the
computational effort to the exact solution
$
X \dvtx[0,T] \times\Omega\rightarrow
{ D ( (-A)^{1/2} ) }
$
of the SPDE~(\ref{reducedForm}) in the sense of
(\ref{particularForm}).
We remark that the log term
in the computational effort $O( N^2\log(N))$ for
$ N \in\mathbb{N} $ arises if one computes the
nonlinearity in~(\ref{reducedForm2}) with fast Fourier
transform (see Figure~\ref{code} for details).

In order to compare the new numerical scheme~(\ref{reducedForm2})
with classical schemes, we consider the well-known
linear implicit Euler scheme
combined with spectral Galerkin methods
applied to the
SPDE~(\ref{reducedForm}). The linear implicit Euler
scheme is denoted by
$\mathcal{F}/\mathcal{B}(H_N)$-measurable mappings
$Z_n^N \dvtx\Omega\rightarrow H_N$,
$n \in\{0, 1, \ldots, N^2 \}$, $N \in\mathbb{N}$, given by
$Z_0^N (\omega) := P_N(\xi(\omega)) + P_N(O_0(\omega))$
and
%
%
\begin{eqnarray} \label{implEulerScheme}\qquad
Z_{n+1}^N (\omega) &:=& \biggl( I - \frac{T}{N^2} A \biggr)^{-1} \biggl\{ Z_n^N (\omega) +
\frac{T}{N^2} (P_N F) ( Z_n^N (\omega) )
\nonumber\\[-8pt]\\[-8pt]
&&\hspace*{73pt}{} + P_N \bigl( B \bigl( W_{(n+1)T/N^2}(\omega) - W_{nT/N^2}(\omega) \bigr) \bigr)
\biggr\} \nonumber
\end{eqnarray}
for every $n \in\{0, 1, \ldots, N^2 - 1 \}$ and
every $N \in\mathbb{N}$.
It has been shown in the literature
(see, e.g.,
Walsh~\cite{w05b},
Gy{\"o}ngy~\cite{g99} and
Hausenblas~\cite{h03a}) that the linear
implicit Euler scheme~(\ref{implEulerScheme})
and other classical numerical schemes
such as the linear implicit Crank--Nicolson scheme
combined with finite elements, finite differences
and spectral Galerkin methods converge with
order $\frac{1}{3}^-$ with respect to the
computational effort.

The following two numerical examples
illustrate the convergence order $\frac{1}{2}^-$ of the
numerical scheme~(\ref{reducedForm2})
and the convergence order $\frac{1}{3}^-$ of the
linear implicit Euler scheme~(\ref{implEulerScheme}).

%
%
%
\begin{figure}[t]

\includegraphics{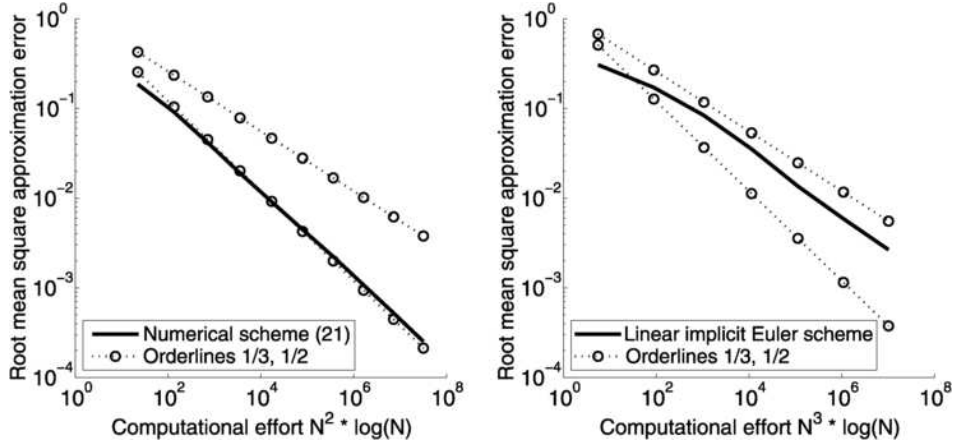}

\caption{Root mean square approximation error (\protect\ref{err1})
of the numerical
scheme (\protect\ref{reducedForm2}) and
root mean square approximation error (\protect\ref{err2})
of the linear
implicit Euler scheme (\protect\ref{implEulerScheme})
applied to SPDE (\protect\ref{exampleformula1}) versus
up to a constant the computational effort.} \label{fig1}
\end{figure}
%
%
%
%
\begin{table}[t]
\caption{Root mean square\vspace*{1pt} approximation error (\protect\ref{err1}) of
$Y_N^{N,N}$ given by (\protect\ref{reducedForm2}) applied~to~the~SPDE~(\protect\ref{exampleformula1}) for $N \in\{ 2^2, 2^3, \ldots, 2^{11}
\}$}\label{tab1}
\begin{tabular*}{\tablewidth}{@{\extracolsep{\fill}}ld{7.0}d{8.0}c@{}}
\hline
& \multicolumn{1}{c}{\textbf{Independent standard}} & \multicolumn{1}{c}{\textbf{Computational}} & \multicolumn{1}{c@{}}{\textbf{Root mean square}}\\
\multicolumn{1}{@{}l}{\textbf{Numerical}} & \multicolumn{1}{c}{\textbf{normal random}} & \multicolumn{1}{c}{\textbf{effort} $\bolds{N^2 \log(N) }$} &
\multicolumn{1}{c@{}}{\textbf{approximation}}\\
\multicolumn{1}{@{}l}{\textbf{scheme~(\ref{reducedForm2})}} & \multicolumn{1}{c}{\textbf{variables} $\bolds{N^2}$} & \multicolumn{1}{c}{\textbf{(up to a constant)}}
& \multicolumn{1}{c@{}}{\textbf{error~(\ref{err1})}}\\
\hline
$Y_{2^2}^{2^2, 2^2}$ & 16 & 22 & 0.1864 \\[4pt]
$Y_{2^3}^{2^3, 2^3}$ & 64 & 133 & 0.0914\\[4pt]
$Y_{2^4}^{2^4, 2^4}$ & 256 & 710 & 0.0417\\[4pt]
$Y_{2^5}^{2^5, 2^5}$ & 1024 & 3 549 & 0.0191\\[4pt]
$Y_{2^6}^{2^6, 2^6}$ & 4096 & 17\mbox{,} 035 & 0.0091\\[4pt]
$Y_{2^7}^{2^7, 2^7}$ & 16\mbox{,} 384 & 79\mbox{,} 496 & 0.0045 \\[4pt]
$Y_{2^8}^{2^8, 2^8}$ & \mathbf{65}\bolds{,}\mathbf{536} & \mathbf{363}\bolds{,}\mathbf{408} & \textbf{0.0022}\\[4pt]
$Y_{2^9}^{2^9, 2^9}$ & 262\mbox{,}144 & 1\mbox{,}635\mbox{,}339 & 0.0011\\[4pt]
$Y_{2^{10}}^{2^{10}, 2^{10}}$ & 1\mbox{,}048\mbox{,}576 & 7\mbox{,}268\mbox{,}174 & 0.0005\\[4pt]
$Y_{2^{11}}^{2^{11}, 2^{11}}$ & 4\mbox{,}194\mbox{,}304 & 31\mbox{,}979\mbox{,}969 & 0.0003\\
\hline
\end{tabular*}
\end{table}

\subsection{A stochastic
reaction diffusion equation} \label{example1}

In this example, we set $\kappa= \frac{1}{100}$,
$T = 1$, $b_n = \frac{n^{-0.55}}{3.5}$ for all $n \in\mathbb{N}$
and consider $f\dvtx[0,1] \times\mathbb{R} \rightarrow\mathbb{R}$
given by $ f(x,y) = 5 \frac{(1-y)}{ (1+y^2) } $ for all $x \in[0,1]$,
$y \in\mathbb{R}$. The SPDE~(\ref{reducedForm}) then reduces to
%
%
\begin{eqnarray} \label{exampleformula1}
dX_t &=& \biggl[ \frac{1}{100} \,\frac{\partial^2}{\partial x^2} X_t + 5 \frac{
(1 - X_t) } { ( 1 + X_t^2 ) } \biggr] \,dt + B \,dW_t,
\nonumber\\
X_t(0) &=& X_t(1) = 0,\\
X_0(x) &=& \frac{ \sin(\pi x) }{ \sqrt{2}} + \frac{3
\sqrt{2}}{5} \sin( 3\pi x ) \nonumber
\end{eqnarray}
for $ x \in[0, 1] $ and $ t \in[0, 1] $.
In Figure~\ref{fig1} (see also
Tables~\ref{tab1} and~\ref{tab2}),
we plot the root mean square discretization
error
%
%
\begin{equation}\label{err1}
\biggl( \mathbb{E} \biggl[ \int_0^1 | X_T (x) - Y_N^{N,N} (x) |^2 \,dx \biggr]
\biggr)^{1/2}
\end{equation}
of the numerical scheme~(\ref{reducedForm2})
versus $ N^2 \log(N) $ (up to a constant
the computational effort)
and the root mean square
discretization error
%
%
\begin{equation}\label{err2}
\biggl( \mathbb{E} \biggl[ \int_0^1 | X_T (x) - Z_{N^2}^{N} (x) |^2 \,dx \biggr]
\biggr)^{1/2}
\end{equation}
of the linear implicit Euler
scheme~(\ref{implEulerScheme})
versus $ N^3 \log(N) $
(up to a constant the computational effort)
for different $ N \in\mathbb{N} $.
The ``expectations'' are
based on $40$ independent random realizations and
the unknown ``exact'' solution is approximated
with a very high accuracy there.
%

%

%
\begin{table}
\caption{Root mean square approximation error (\protect\ref{err2})
of $Z_{N^2}^{N}$ given by (\protect\ref{implEulerScheme})
applied to~the~SPDE~(\protect\ref{exampleformula1}) for
$N \in\{ 2^1, 2^2, \ldots, 2^7 \}$}\label{tab2}
\begin{tabular*}{\tablewidth}{@{\extracolsep{\fill}}ld{7.0}d{8.0}c@{}}
\hline
\textbf{Linear implicit} & \multicolumn{1}{c}{\textbf{Independent standard}} & \multicolumn{1}{c}{\textbf{Computational}} & \multicolumn{1}{c@{}}{\textbf{Root mean square}}\\
\textbf{Euler} & \multicolumn{1}{c}{\textbf{normal random}} & \multicolumn{1}{c}{\textbf{effort} $\bolds{N^3\log(N)}$} & \multicolumn{1}{c@{}}{\textbf{approximation}}\\
\textbf{scheme~(\ref{implEulerScheme})} & \multicolumn{1}{c}{\textbf{variables} $\bolds{N^3}$} & \multicolumn{1}{c}{\textbf{(up to a constant)}} & \multicolumn{1}{c@{}}{\textbf{error
(\ref{err2})}}\\
\hline
$Z_{2^2}^{2^1}$ & 8 & 6 & 0.3066\\[4pt]
$Z_{2^4}^{2^2}$ & 64 & 88 & 0.1715\\[4pt]
$Z_{2^6}^{2^3}$ & 512 & 1 064 & 0.0837\\[4pt]
$Z_{2^8}^{2^4}$ & 4 096 & 11\mbox{,} 356 & 0.0353\\[4pt]
$Z_{2^{10}}^{2^5}$ & 32\mbox{,}768 & 113\mbox{,} 565 & 0.0135\\[4pt]
$Z_{2^{12}}^{2^6}$ & 262\mbox{,}144 & 1\mbox{,}090\mbox{,}226 & 0.0058\\[4pt]
$Z_{2^{14}}^{2^7}$ & \mathbf{2}\bolds{,}\mathbf{097}\bolds{,}\mathbf{152} & \mathbf{10}\bolds{,}\mathbf{175}\bolds{,}\mathbf{444} & \textbf{0.0027}\\
\hline
\end{tabular*}
\end{table}

The short \textsc{matlab} code in Figure~\ref{code} shows that
the solution of SPDE~(\ref{exampleformula1})
can be simulated quite easily with the
numerical scheme~(\ref{reducedForm2}).
Figure~\ref{fig1verbatim} is the result of the \textsc{matlab}
code in Figure~\ref{code}.
It shows the solution of the stochastic
reaction diffusion equation~(\ref{exampleformula1})
at time $t = T = 1$ for one sample path $\omega\in\Omega$
approximated with the numerical
method~(\ref{reducedForm2}).

%
%
\begin{figure}[b]
{\fontsize{8pt}{10pt}\selectfont{
\begin{verbatim}
N = 1000; T = 1; A = - pi^2 * (1:N).^2 / 100; Y = [1/2,0,3/5,zeros(1,N-3)];
S = sqrt( ( exp(2*T/N*A) - 1 ) ./ A / 2 ) / 3.5 .* (1:N).^ -0.55;
for n=1:N
  y = dst(Y) * sqrt(2);
  FY = idst( 5 * ( 1 - y ) ./ ( 1 + y.^2 ) ) / sqrt(2);
  Y = exp( A * T/N ) .* ( Y + T/N * FY ) + S .* randn(1,N);
end
plot( (0:N+1)/(N+1), [0,dst(Y)*sqrt(2),0], 'k', 'Linewidth', 2 );
\end{verbatim}}}
\caption{\textsc{matlab} code for the numerical scheme
(\protect\ref{reducedForm2}) applied to the SPDE
(\protect\ref{exampleformula1}).}\label{code}\vspace*{-4pt}
\end{figure}

%
%
\begin{figure}

\includegraphics{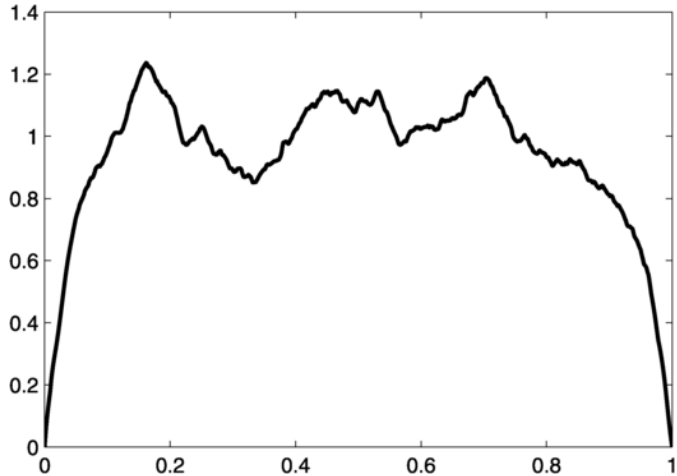}

\caption{Result of the \textsc{matlab}
code in Figure \protect\ref{code}:
Solution of the stochastic reaction
diffusion equation (\protect\ref{exampleformula1}) at $t = T = 1$
for one sample path $\omega\in\Omega$ approximated with the
numerical
method (\protect\ref{reducedForm2}).}\label{fig1verbatim}\vspace*{-2pt}
\end{figure}

%

\subsection{A stochastic partial differential
equation with a spatially dependent $f$}

\begin{figure}

\includegraphics{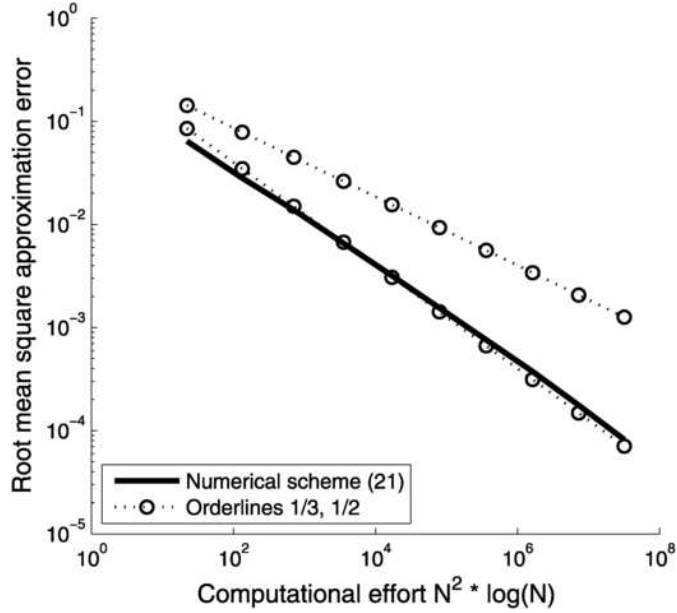}

\caption{Root mean square approximation error (\protect\ref{err1})
of the numerical scheme (\protect\ref{reducedForm2})
applied to SPDE (\protect\ref{exampleformula3}) versus
up to a constant the computational effort.}\label{fig3}
\end{figure}

\begin{figure}

\includegraphics{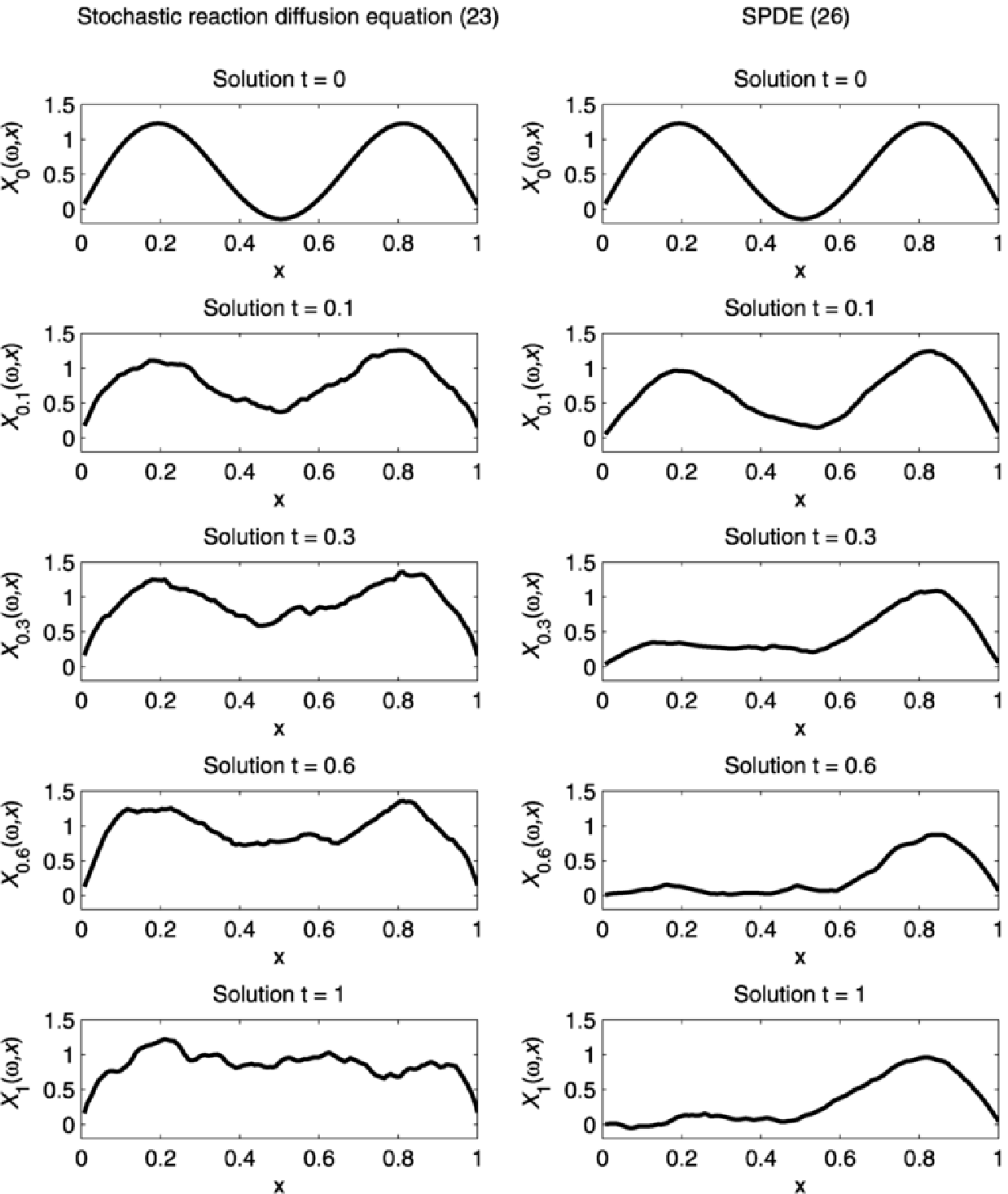}

\caption{Solution $X_t( \omega, x )$, $x \in[0,1]$,
of the stochastic reaction diffusion equation
(\protect\ref{exampleformula1})
and of the SPDE (\protect\ref{exampleformula3})
for $t \in\{0, \frac{1}{10}, \frac{3}{10},
\frac{6}{10},1 \}$ and one sample path
$\omega\in\Omega$ approximated with the numerical
method (\protect\ref{reducedForm2}).}\label{subplots}
\end{figure}

This time let $\kappa= \frac{1}{50}$, $T = 1$,
$b_n = \frac{n^{-0.6}}{5}$ for all $n \in\mathbb{N}$
and consider
$f\dvtx[0,1] \times\mathbb{R} \rightarrow\mathbb{R}$
given by
$ f(x,y) = (3.8 x^2 -2) y $ for all $x \in[0,1]$,
$y \in\mathbb{R}$ to obtain the SPDE
%
%
\begin{eqnarray} \label{exampleformula3}
dX_t &=& \biggl[ \frac{1}{50} \,\frac{\partial^2}{\partial x^2} X_t + (3.8 x^2 -
2) X_t \biggr] \,dt + B \,dW_t,
\nonumber\\
X_t(0) &=& X_t(1) = 0,\\
X_0(x) &=& \frac{ \sin(\pi x) }{ \sqrt{2}} - \frac{3
\sqrt{2}}{5} \sin( 3\pi x ) \nonumber
\end{eqnarray}
for $ x \in[0, 1] $ and $ t \in[0, T] $.
Here too,\vspace*{-1pt} the numerical approximation
(\ref{reducedForm2}) converges to the exact solution with\vspace*{1pt}
order $\frac{1}{2}^-$ with respect to up to a constant the
computational effort (see Figure~\ref{fig3}).
Finally, in Figure~\ref{subplots} we illustrate how the two different
$f$ from examples~(\ref{exampleformula1}) and~(\ref{exampleformula3})
affect the evolution of the respective solution $X_t(\omega,x)$, $ x
\in[0,1] $, for $ t \in\{ 0, \frac{1}{10}, \frac{3}{10},
\frac{6}{10}, 1 \} $ and one sample path
$\omega\in\Omega$.

%
%
%
%

%
%

\section{A further numerical example} \label{sec.furtherexamples}

\begin{figure}

\includegraphics{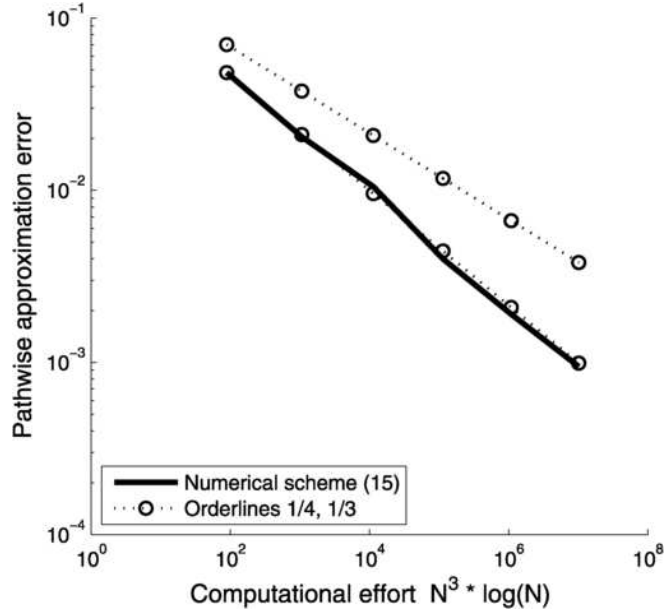}

\caption{Pathwise approximation error
of the numerical scheme (\protect\ref{expEulerScheme})
applied to SPDE~(\protect\ref{2dformula}) versus
up to a constant the computational effort for one random
$\omega\in\Omega$.}\label{fig_error_2d}
\end{figure}

Although our setting in Section~\ref{sec.assump} uses
the standard global Lipschitz assumption on the nonlinearity
of the SPDE, we demonstrate the efficiency of our method numerically
for a SPDE with a cubic nonglobally Lipschitz nonlinearity in
this section.
More formally, we consider the SPDE
%
%
\begin{equation} \label{2dformula}
dX_t = \biggl[ \frac{1}{10} \biggl( \frac{\partial^2}{\partial x_1^2} +
\frac{\partial^2}{\partial x_2^2} \biggr) X_t + X_t - X_t^3 \biggr] \,dt + dW^Q_t,
\end{equation}
with
\[
X_t |_{\partial(0,1)^2} \equiv0
\]
and
\[
X_0(x_1, x_2) = \sin(\pi x_1)
\sin( \pi x_2 )
\]
for $ t, x_1, x_2 \in[0, 1] $ on the $\mathbb{R}$-Hilbert space
$H = L^2 ( (0,1)^2, \mathbb{R} ) $
of equivalence classes of
$\mathcal{B} ((0,1)^2)/\mathcal{B}(\mathbb{R})$-measurable
and square integrable functions from $(0,1)^2$ to $\mathbb{R}$
here where $( W^Q_t )_{t \in[0,1]} $ is a cylindrical
$Q$-Wiener process on $H$ with the covariance operator
$Q\dvtx H \rightarrow H$ given by
\begin{eqnarray*}
&& (Qv)(x_1, x_2)
\\
&&\qquad
= \sum_{n,m = 1}^\infty\frac{4 \sin(n \pi x_1) \sin(n \pi x_2)
}{(n+m)^2} \\
&&\qquad\quad\hspace*{23.2pt}{}\times\int_0^1 \int_0^1 \sin(n\pi y_1) \sin(m \pi y_2) v(y_1, y_2)
\,dy_1 \,dy_2
\end{eqnarray*}
for all $x_1, x_2 \in(0,1)$ and all $v \in H$.
Of course,~(\ref{2dformula}) is not included in our
setting in Section~\ref{sec.assump}. Even worse, it has recently
been shown in~\cite{hj09a} that many numerical methods
fail to converge to the solution of a stochastic
differential equation with super linearly growing coefficients
in the strong root mean square sense. However, convergence in the
pathwise sense often holds due to Gy{\"o}ngy's result~\cite{g98b}.
Therefore,
we plot in Figure~\ref{fig_error_2d} the pathwise difference
\[
\biggl( \int_0^1 \int_0^1 | X_T(\omega,x_1,x_2) -
Y_N^{N,N} (\omega, x_1, x_2) |^2 \,dx_1 \,dx_2 \biggr)^{1/2}
\]
of the exact solution $X_T(\omega)$ and of the numerical
approximation $Y_N^{N,N} (\omega)$ [see~(\ref{expEulerScheme})]
applied to the SPDE~(\ref{2dformula}) versus up to a constant
the computational effort $N^3 \log(N)$ for $N \in\{ 2^2, 2^3, \ldots,
2^7 \}$ and
one random $\omega\in\Omega$. It turns out that the
method~(\ref{expEulerScheme}) converges with order $\frac{1}{3}^-$
with respect
to the computational effort. The linear implicit Euler scheme is known
to converge in the pathwise sense with order $\frac{1}{4}^-$ with
respect to
%
%
the computational effort to the solution of the SPDE (\ref
{2dformula}). Further
pathwise approximation results for the SPDE~(\ref{2dformula}) and
other SPDEs with nonglobally Lipschitz coefficients can be found in
\cite{g99,gm09,gm05} and~\cite{ps05},
for instance.
%
%
\begin{figure}

\includegraphics{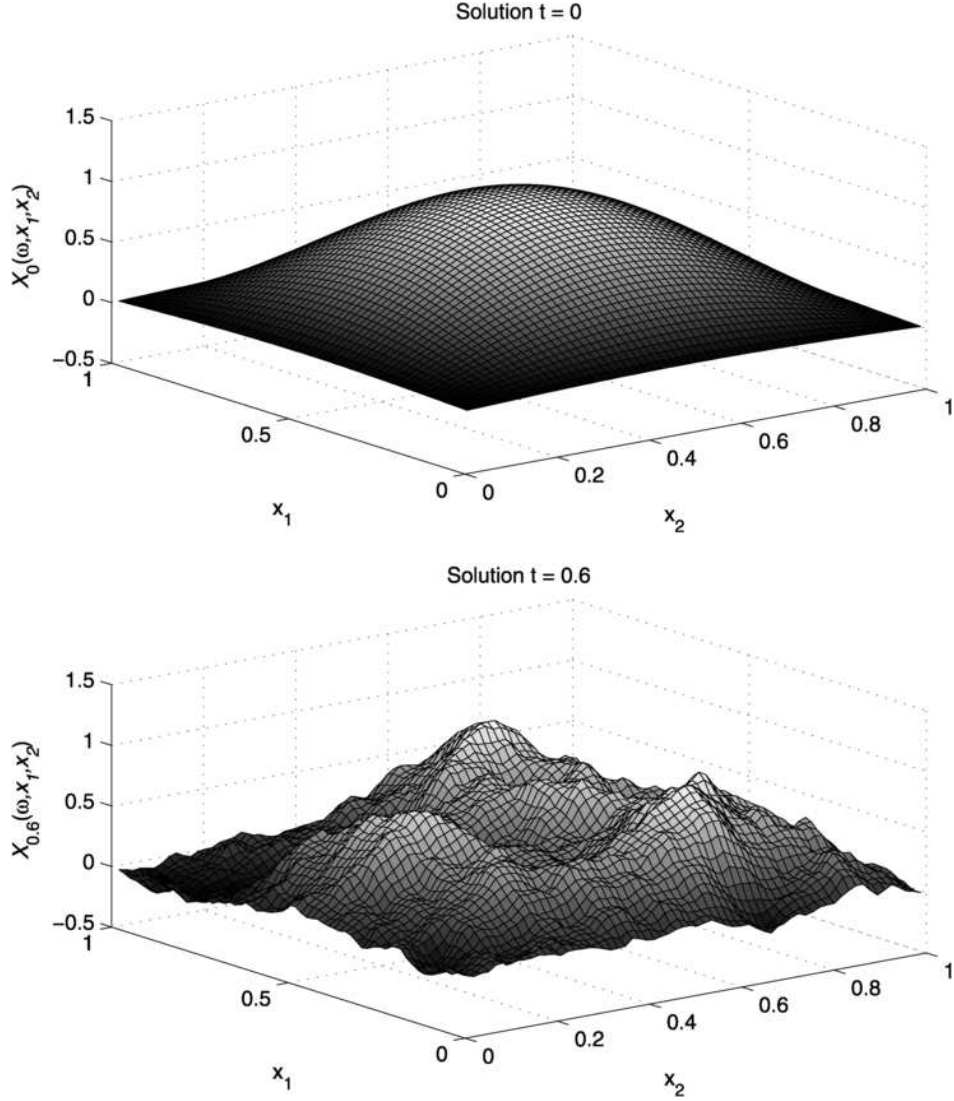}

\caption{Solution $X_t( \omega, x_1, x_2 )$, $x_1, x_2 \in[0,1]$,
of the SPDE (\protect\ref{2dformula})
for $t \in\{0, \frac{6}{10} \}$ and one random
$\omega\in\Omega$ approximated with the numerical
method (\protect\ref{expEulerScheme}).}\label{fig_2d_surface}
\end{figure}
Finally, we plot the solution of SPDE~(\ref{2dformula})
for $t \in\{0, \frac{6}{10}\}$
and one random $\omega\in\Omega$ in Figure~\ref{fig_2d_surface}.

\section{Proofs} \label{proofs}
The notation
\[
\| Z \|_{L^p (\Omega; W)} := ( \mathbb{E} [ \| Z \|_W^p ] )^{1/p}
\ \in[0,\infty]
\]
is used throughout this section for an $ \mathbb{R} $-Banach
space $ ( W, \| \cdot\|_W ) $, a $ \mathcal{F}
/ \mathcal{B}(W) $-measurable
mapping $ Z \dvtx\Omega\rightarrow W $ and a real
number $ p \in[1,\infty) $.
%
\subsection{\texorpdfstring{Proof of Theorem \protect\ref{theorem1}}{Proof of Theorem 1}}
The
$ \mathcal{F} / \mathcal{B}(H) $-measurable
mappings $ Y^M_m \dvtx\Omega\rightarrow H $ for
$ m \in\{ 0, 1, \ldots, M \} $ and
$ M \in\mathbb{N} $ given by
%
%
\begin{eqnarray} \label{NumScheme1}
Y^M_m(\omega) &:=& e^{A m h } \xi(\omega) + h \Biggl( \sum^{m-1}_{k=0} e^{A ( m
h - k h ) } F( X_{ k h }(\omega) ) \Biggr) \nonumber\\[-8pt]\\[-8pt]
&&{}+ O_{ m h }(\omega)\nonumber
\end{eqnarray}
for every $ m \in\{ 0, 1, \ldots, M \} $, $ \omega\in\Omega$ and $ M
\in\mathbb{N} $ are used throughout this proof. Here and below $ h $ is
the time stepsize $ h = h_M = \frac{T} {M} $ with $ M \in\mathbb{N} $.
This proof is divided\vadjust{\goodbreak} into three parts. In the first part (see Section
\ref{temp_error}), we estimate
\[
\| X_{m h} - Y_m^M \|_{ L^2 ( \Omega; H ) }
\]
for every $ m \in\{ 0, 1, \ldots, M \} $ and every $ M \in\mathbb{N} $
which corresponds to the temporal discretization error. In the second
part (see Section~\ref{spat_error}), we estimate
\[
\| Y_m^M - P_N( Y_m^M ) \|_{ L^2 ( \Omega; H ) }
\]
for every $ m \in\{0, 1, \ldots, M \} $ and every $ N, M \in\mathbb{N}$
which corresponds to the spatial discretization error. Finally, we
estimate
\[
\| P_N( Y_m^M ) - Y_m^{N, M} \|_{ L^2 ( \Omega; H ) }
\]
for every $ m \in\{ 0, 1, \ldots, M \} $ and every $ N, M \in\mathbb{N}
$ in the third part (see Section~\ref{refLipschitzEstimate}). Combining
these three parts will then yield the desired assertion via Gronwall's lemma as
we will see below.

Before we begin with the first part, we introduce a universal constant
$R > 0$ which is needed throughout this proof. More precisely, let $ R
\in(0, \infty) $ be a real number which satisfies
\begin{eqnarray*}
\| F( X_t ) \|_{ L^2 ( \Omega; H ) } &\leq& R,\\
\| ( X_{t_2} - O_{t_2} )-
( X_{t_1} - O_{t_1} ) \|_{ L^2 ( \Omega; H ) } &\leq& R | t_2 - t_1 | ,
\\
\frac{1}{\lambda_1} + \frac{1}{(1-\gamma)} + T + c &\leq& R,\\
\| v \|_H &\leq& R \| v \|_V,\\
\| O_{t_2} - O_{t_1} \|_{ L^4 ( \Omega; V ) } &\leq& R | t_2 - t_1 |^{
\theta} ,
\\
\| \xi\|_{ L^2 ( \Omega; D((-A)^\gamma) } &\leq& R,\\
\| O_t \|_{ L^4 (
\Omega; D( (-A)^\gamma) ) } &\leq& R,\\
\| X_t \|_{ L^4( \Omega; D(
(-A)^{ {1/2} } ) ) } &\leq& R
\end{eqnarray*}
for every $t, t_1, t_2 \in[0, T ] $ and every $ v \in V $ where $
\lambda_1 \in(0,\infty) $ is given in Assumption~\ref{semigroup} where
$ c \in[0,\infty)$ is given in Assumption~\ref{drift} and where
$\gamma\in[ \frac{1}{2},1 )$ and $\theta\in( 0, \frac{1}{2} ] $ are
given in Assumption~\ref{stochconv}. Indeed, such a real number exists
due to Assumptions~\ref{semigroup}--\ref{initial} and Lemma
\ref{lemmaRef3} in Section~\ref{secprop}.
%
\subsubsection{Temporal discretization error}\label{temp_error}
Due to~(\ref{eqSolution}), we have
\begin{eqnarray*}
X_{m h} &=& e^{ A {m h} } \xi+ \int_0^{m h} e^{ A {( m h - s ) } } F(
X_s ) \,ds + O_{m h}
\\
&=&
e^{ A {mh} } \xi
+ \sum_{ k = 0 }^{ m - 1 }
\int_{ k h }^{( k + 1 ) h}
e^{ A {( m h - s ) } }
F( X_s ) \,ds
+ O_{mh}
\end{eqnarray*}
for every $ m \in\{0, 1, \ldots, M\} $ and
every $ M \in\mathbb{N} $.
From~(\ref{NumScheme1}), we have
\begin{eqnarray*}
&& \| X_{m h} - Y_m^M \|_{L^2 ( \Omega; H ) }
\\
&&\qquad
= \Biggl\| \sum_{k = 0}^{m - 1} \int_{k h}^{( k + 1 ) h } e^{ A {( m h - s )}
} F(X_s) \,ds - h \Biggl( \sum_{k = 0}^{m - 1} e^{ A { ( m h - k h ) } } F(X_{k
h}) \Biggr) \Biggr\|_{L^2 ( \Omega; H ) }
\\
&&\qquad\leq
\Biggl\| \sum_{k = 0}^{m - 2} \int_{k h}^{( k + 1 ) h } e^{ A {( m h - s ) }
} F(X_s) \,ds - h \Biggl( \sum_{k = 0}^{m - 2} e^{ A { ( m h - k h ) } } F(X_{k
h}) \Biggr) \Biggr\|_{L^2 ( \Omega; H ) }
\\[-2pt]
&&\qquad\quad{}
+ \int_{ \max ( m - 1, 0 ) h  }^{m h} \bigl\| e^{ A {( m h - s ) } }
F(X_s) \bigr\|_{L^2 ( \Omega; H ) } \,ds
\\[-2pt]
&&\qquad\quad{}
+ h \bigl\| e^{ A h } F\bigl( X_{ \max ( m - 1, 0 ) h  } \bigr) \bigr\|_{L^2 ( \Omega; H
) }
\end{eqnarray*}
and
\begin{eqnarray*}
&& \| X_{m h} - Y_m^M \|_{L^2 ( \Omega; H ) }
\\[-2pt]
&&\qquad
\leq\Biggl\| \sum_{k = 0}^{m - 2} \int_{k h}^{( k + 1 ) h } e^{ A {( m h - s
) } } F(X_s) \,ds - h \Biggl( \sum_{k = 0}^{m - 2} e^{ A { ( m h - k h ) } }
F(X_{k h}) \Biggr) \Biggr\|_{L^2 ( \Omega; H ) }
\\[-2pt]
&&\qquad\quad{}
+ \int_{ \max( m - 1, 0 ) h }^{m h} \bigl\| e^{ A {( m h - s ) } } \bigr\|_{L (
H ) } \| F(X_s) \|_{ L^2 ( \Omega; H ) } \,ds
\\[-2pt]
&&\qquad\quad{}
+ h \| e^{ A h } \|_{ L ( H ) } \bigl\| F\bigl( X_{ \max ( m - 1, 0 ) h } \bigr)
\bigr\|_{ L^2 ( \Omega; H ) }
\end{eqnarray*}
for every $ m \in\{ 0,1,\ldots, M \}$
and every $M \in\mathbb{N} $.
Therefore, we obtain
\begin{eqnarray*}
&& \| X_{m h} - Y_m^M \|_{ L^2 ( \Omega; H ) }
\\[-2pt]
&&\qquad\leq
\Biggl\| \sum_{k = 0}^{m - 2} \int_{k h}^{( k + 1 ) h } e^{ A {( m h - s ) }
} F(X_s) \,ds\\[-2pt]
&&\qquad\quad\hspace*{29.1pt}\hspace*{-18.2pt}{} - h \Biggl( \sum_{k = 0}^{m - 2} e^{ A { ( m h - k h ) } } F(X_{k
h}) \Biggr) \Biggr\|_{L^2 ( \Omega; H ) }
+ 2 R h
\\[-2pt]
&&\qquad\leq
\Biggl\| \sum_{k = 0}^{m - 2} \int_{k h}^{( k + 1 ) h} e^{ A {( m h - s )} }
\bigl( F(X_s) - F(X_{k h}) \bigr) \,ds \Biggr\|_{L^2 ( \Omega; H ) }\\[-2pt]
&&\qquad\quad{}
 + 2 R h
\\[-2pt]
&&\qquad\quad{}
+ \Biggl\| \sum_{k = 0}^{m - 2} \int_{k h}^{( k + 1 ) h } e^{ A {( m h - s )}
} F(X_{k h}) \,ds \\[-2pt]
&&\qquad\quad\hspace*{29.1pt}{} - h \Biggl( \sum_{k = 0}^{m - 2} e^{ A { ( m h - k h ) } }
F(X_{k h}) \Biggr) \Biggr\|_{L^2 ( \Omega; H ) }
\end{eqnarray*}
and
\begin{eqnarray*}
&& \| X_{m h} - Y_m^M \|_{ L^2 ( \Omega; H ) }
\\[-2pt]
&&\qquad\leq
\Biggl\| \sum_{k = 0}^{m - 2} \int_{k h}^{( k + 1 ) h} e^{ A {( m h - s )} }
\bigl( F(X_s) - F(X_{k h} + O_s - O_{k h}) \bigr) \,ds \Biggr\|_{ L^2 ( \Omega; H ) }
\\[-2pt]
&&\qquad\quad{}
+ \Biggl\| \sum_{k = 0}^{m - 2} \int_{k h}^{( k + 1 ) h } e^{ A {( m h - s )}
} \bigl( F(X_{k h} + O_s - O_{k h}) - F(X_{k h}) \bigr) \,ds \Biggr\|_{ L^2 ( \Omega; H )
}
\\[-2pt]
&&\qquad\quad{}
+ \sum_{k = 0}^{m - 2} \biggl\| \int_{kh}^{(k + 1 ) h} \bigl( e^{ A {( m h - s )}
} - e^{ A {( m h - k h )} } \bigr) F(X_{kh}) \,ds \biggr\|_{L^2 ( \Omega; H ) }
+ 2R h
\end{eqnarray*}
for every $ m \in\{ 0,1,\ldots, M \} $
and every $ M \in\mathbb{N} $. Hence, we obtain
\begin{eqnarray*}
&& \| X_{m h} - Y_m^M \|_{L^2 ( \Omega; H ) }
\\[-0.5pt]
&&\qquad\leq
\sum_{k = 0}^{m - 2} \int_{k h}^{( k + 1 ) h} \bigl\| e^{ A {( m h - s )} }
\bigr\|_{ L( H ) } \| F(X_s) - F(X_{k h} + O_s - O_{k h}) \|_{ L^2 ( \Omega;
H ) } \,ds
\\[-0.5pt]
&&\qquad\quad{}
+ \Biggl\| \sum_{k = 0}^{m - 2} \int_{k h}^{( k + 1 ) h} e^{ A { ( m h - s )
} } F'( X_{ k h } ) ( O_s - O_{ k h } ) \,ds \Biggr\|_{ L^2 ( \Omega; H ) }
\\[-0.5pt]
&&\qquad\quad{}
+ \Biggl\| \sum_{k = 0}^{m - 2} \int_{k h}^{( k + 1 ) h } e^{ A {( m h - s )}
} \int_0^1 F'' \bigl( X_{k h} + r ( O_s - O_{k h}) \bigr)
\\[-0.5pt]
&&\qquad\quad\hspace*{143.3pt}\hspace*{-15.6pt}{}
\times( O_s - O_{ k h }, O_s - O_{ k h }
) \\[-0.5pt]
&&\qquad\quad\hspace*{172pt}{}\times( 1 - r ) \,dr \,ds \Biggr\|_{ L^2 ( \Omega; H ) }
\\[-0.5pt]
&&\qquad\quad{}
+ \sum_{k = 0}^{m - 2} \int_{k h}^{( k + 1 ) h } \bigl\| \bigl( e^{ A {( m h - s
)} } - e^{ A {( m h - k h )} } \bigr) F(X_{k h}) \bigr\|_{ L^2 ( \Omega; H ) } \,ds
+ 2 R^2 M^{-1}
\end{eqnarray*}
and
\begin{eqnarray*}
&& \| X_{m h} - Y_m^M \|_{L^2 ( \Omega; H ) }
\\[-0.5pt]
&&\qquad\leq
c \sum_{k = 0}^{m - 2} \int_{k h}^{( k + 1 ) h} \| X_s - (X_{k h} + O_s
- O_{k h}) \|_{ L^2 ( \Omega; H ) } \,ds
\\[-0.5pt]
&&\qquad\quad{}
+ \Biggl\| \sum_{k = 0}^{m - 2} \int_{k h}^{( k + 1 ) h} e^{ A { ( m h - s )
} } F'( X_{ k h } ) \bigl( O_s - e^{ A { ( s - k h ) } } O_{ k h } \bigr) \,ds \Biggr\|_{
L^2 ( \Omega; H ) }
\\[-0.5pt]
&&\qquad\quad{}
+ \Biggl\| \sum_{k = 0}^{m - 2} \int_{k h}^{( k + 1 ) h} e^{ A { ( m h - s )
} } F'( X_{ k h } ) \bigl( \bigl( e^{ A { ( s - k h ) } } - I \bigr) O_{ k h } \bigr) \,ds
\Biggr\|_{ L^2 ( \Omega; H ) }
\\[-0.5pt]
&&\qquad\quad{}
+ \sum_{k = 0}^{m - 2} \int_{k h}^{( k + 1 ) h } \bigl\| e^{ A {( m h - s )}
} \bigr\|_{ L ( H ) } \\[-0.5pt]
&&\qquad\quad\hspace*{65.7pt}{}\times\int_0^1 \bigl\| F'' \bigl( X_{k h} + r ( O_s - O_{k h})
\bigr)\\[-0.5pt]
&&\qquad\quad\hspace*{96.5pt}{}\times
( O_s - O_{ k h }, O_s - O_{ k h } ) \bigr\|_{ L^2 ( \Omega; H ) }
\,dr \,ds
\\[-0.5pt]
&&\qquad\quad{}
+ \sum_{k = 0}^{m - 2} \int_{k h}^{( k + 1 ) h } \bigl\| e^{ A {( m h - s )}
} - e^{ A {( m h - k h )} } \bigr\|_{ L ( H ) } \| F(X_{k h}) \|_{ L^2 (
\Omega; H ) } \,ds
\\[-0.5pt]
&&\qquad\quad{}+ 2 R^2 M^{-1}
\end{eqnarray*}
for every $ m \in\{ 0,1,\ldots, M \} $
and every $ M \in\mathbb{N} $. Therefore, we have
%
%
\begin{eqnarray} \label{eqPart5}
&& \| X_{m h} - Y_m^M \|_{ L^2 ( \Omega; H ) } \nonumber
\\
&&\qquad \leq
c \sum_{k = 0}^{m - 2} \int_{k h}^{( k + 1 ) h} \| ( X_s - O_s ) - (
X_{k h} - O_{k h} ) \|_{ L^2 ( \Omega; H ) } \,ds \nonumber
\\
&&\qquad\quad{}
+ \Biggl\| \sum_{k = 0}^{m - 2} \int_{k h}^{( k + 1 ) h} e^{ A { ( m h - s )
} } F'(X_{k h}) \bigl( O_s - e^{ A { ( s - kh ) } } O_{k h} \bigr) \,ds \Biggr\|_{ L^2 (
\Omega; H ) } \nonumber\\
&&\qquad\quad{}
+ \sum_{k = 0}^{m - 2} \int_{k h}^{( k + 1 ) h} \bigl\| e^{ A { ( m h - s )
} } F'( X_{ k h } ) \bigl( \bigl( e^{ A { ( s - k h ) } } - I \bigr) O_{ k h } \bigr) \bigr\|_{
L^2 ( \Omega; H ) } \,ds
\\
&&\qquad\quad{}
+ c R \sum_{k = 0}^{m - 2} \int_{k h}^{( k + 1 ) h } \int_0^1 \bigl\| \| O_s
- O_{k h} \|_V^2 \bigr\|_{ L^2 ( \Omega; \mathbb{R} ) } \,dr \,ds \nonumber
\\
&&\qquad\quad{}
+ \sum_{k = 0}^{m - 2} \int_{k h}^{( k + 1 ) h } \frac{ ( m h - k h - m
h + s ) }{ ( m h - s ) } \| F(X_{k h}) \|_{ L^2 ( \Omega; H ) } \,ds\nonumber\\
&&\qquad\quad{} + 2
R^2 M^{-1} \nonumber
\end{eqnarray}
for every $ m \in\{ 0,1,\ldots, M \} $
and every $ M \in\mathbb{N} $ due to
Lemma~\ref{lemmaRef4} below (see Section~\ref{secprop}).
Furthermore, we have
\begin{eqnarray*}
&&\mathbb{E} \Biggl[ \Biggl\| \sum_{k = 0}^{m - 2} \int_{k h}^{( k + 1 ) h} e^{ A { (
m h - s ) } } F'(X_{k h}) \bigl( O_s - e^{ A { ( s - kh ) } } O_{k h} \bigr) \,ds
\Biggr\|_H^2 \Biggr]
\\ 
&&\qquad= \sum_{k,\tilde{k} = 0}^{m - 2} \mathbb{E} \biggl[ \biggl\langle \int_{k h}^{( k + 1 )
h} e^{ A { ( m h - s ) } } F'(X_{k h}) \bigl( O_s - e^{ A { ( s - kh ) } }
O_{k h} \bigr) \,ds,
\\
&&\qquad\quad\hspace*{43.2pt}\int_{\tilde{k} h}^{( \tilde{k} + 1 ) h} e^{ A { ( m h - s ) } }
F'(X_{\tilde{k} h}) \bigl( O_s - e^{ A { ( s - \tilde{k} h ) } }
O_{\tilde{k} h} \bigr) \,ds
\biggr\rangle_H
\biggr]
\end{eqnarray*}
and hence
\begin{eqnarray*}
&&\mathbb{E} \Biggl[ \Biggl\| \sum_{k = 0}^{m - 2} \int_{k h}^{( k + 1 ) h} e^{ A { (
m h - s ) } } F'(X_{k h}) \bigl( O_s - e^{ A { ( s - kh ) } } O_{k h} \bigr) \,ds
\Biggr\|_H^2 \Biggr]
\\
&&\qquad= \sum_{k=0}^{m-2} \mathbb{E} \biggl[ \biggl\| \int_{k h}^{( k + 1 ) h} e^{ A { ( m
h - s ) } } F'(X_{k h}) \\
&&\qquad\quad\hspace*{71.6pt}{}\times \bigl( O_s - e^{ A { ( s - kh ) } } O_{k h} \bigr) \,ds
\biggr\|_H^2 \biggr]
\\
&&\qquad\quad{} + \mathop{\sum_{k,\tilde{k} = 0}}_{k \neq\tilde{k}}^{m - 2} \mathbb{E}
\biggl[ \biggl\langle\int_{k h}^{( k + 1 ) h} e^{ A { ( m h - s ) } } F'(X_{k h})\bigl( O_s
- e^{ A { ( s - kh ) } } O_{k h} \bigr) \,ds,
\\
&&\qquad\quad\hspace*{56.4pt}\int_{\tilde{k} h}^{( \tilde{k} + 1 ) h} e^{ A { ( m h - s ) } }
F'(X_{\tilde{k} h}) \bigl( O_s - e^{ A { ( s - \tilde{k} h ) } }
O_{\tilde{k} h} \bigr) \,ds
\biggr\rangle_H\biggr]
\end{eqnarray*}
for every $ m \in\{ 0,1,\ldots, M \} $
and every $ M \in\mathbb{N} $.
This yields
\begin{eqnarray*}
&&\mathbb{E} \Biggl[ \Biggl\| \sum_{k = 0}^{m - 2} \int_{k h}^{( k + 1 ) h} e^{ A { (
m h - s ) } } F'(X_{k h}) \bigl( O_s - e^{ A { ( s - kh ) } } O_{k h} \bigr) \,ds
\Biggr\|_H^2 \Biggr]
\\
&&\qquad= \sum_{k=0}^{m-2} \mathbb{E} \biggl[ \biggl\| \int_{k h}^{( k + 1 ) h} e^{ A { ( m
h - s ) } } F'(X_{k h}) \bigl( O_s - e^{ A { ( s - kh ) } } O_{k h} \bigr) \,ds
\biggr\|_H^2 \biggr]
\\
&&\qquad\quad{}+ 2 \mathop{\sum_{k,\tilde{k} = 0}}_{k < \tilde{k}}^{m - 2} \mathbb{E}
\biggl[ \biggl\langle \int_{k h}^{( k + 1 ) h} e^{ A { ( m h - s ) } } F'(X_{k h})\\
&&\qquad\quad\hspace*{76.4pt}\hspace*{18.7pt}{}\times \bigl( O_s
- e^{ A { ( s - kh ) } } O_{k h} \bigr) \,ds,
\\
&&\qquad\quad\hspace*{63.3pt}\hspace*{-3.4pt}\int_{\tilde{k} h}^{( \tilde{k} + 1 ) h} e^{ A { ( m h - s ) } }
F'(X_{\tilde{k} h})
\\
&&\qquad\quad\hspace*{76.4pt}\hspace*{21pt}{}\times \bigl( O_s - e^{ A { ( s - \tilde{k} h ) } }
O_{\tilde{k} h} \bigr) \,ds\biggr\rangle_H\biggr]
\end{eqnarray*}
and
\begin{eqnarray*}
&&\mathbb{E} \Biggl[ \Biggl\| \sum_{k = 0}^{m - 2} \int_{k h}^{( k + 1 ) h} e^{ A { (
m h - s ) } } F'(X_{k h}) \bigl( O_s - e^{ A { ( s - kh ) } } O_{k h} \bigr) \,ds
\Biggr\|_H^2 \Biggr]
\\
&&\qquad= \sum_{k=0}^{m-2} \mathbb{E} \biggl[ \biggl\| \int_{k h}^{( k + 1 ) h} e^{ A { ( m
h - s ) } } F'(X_{k h}) \bigl( O_s - e^{ A { ( s - kh ) } } O_{k h} \bigr) \,ds
\biggr\|_H^2 \biggr]
\\
&&\qquad\quad{}+ 2 \mathop{\sum_{k,\tilde{k} = 0}}_{k < \tilde{k}}^{m - 2} \mathbb{E}
\biggl[ \mathbb{E} \biggl[ \biggl\langle \int_{k h}^{( k + 1 ) h} e^{ A { ( m h - s ) } }
F'(X_{k h}) \\
&&\qquad\quad\hspace*{76.4pt}\hspace*{30.6pt}{}\times\bigl( O_s - e^{ A { ( s - kh ) } } O_{k h} \bigr) \,ds,
\\
&&\qquad\quad\hspace*{76.4pt}\hspace*{-3.1pt}\int_{\tilde{k} h}^{( \tilde{k} + 1 ) h} e^{ A { ( m h - s ) } }
F'(X_{\tilde{k} h}) \\
&&\qquad\quad\hspace*{76.4pt}\hspace*{34pt}{}\times\bigl( O_s - e^{ A { ( s - \tilde{k} h ) } }
O_{\tilde{k} h} \bigr) \,ds
\biggr\rangle_H\Big|\mathcal{F}_{\tilde{k} h}\biggr]\biggr]
\end{eqnarray*}
for every $ m \in\{ 0,1,\ldots, M \} $
and every $ M \in\mathbb{N} $.
Hence, we obtain
\begin{eqnarray*}
&&\mathbb{E} \Biggl[ \Biggl\| \sum_{k = 0}^{m - 2} \int_{k h}^{( k + 1 ) h} e^{ A { (
m h - s ) } } F'(X_{k h}) \bigl( O_s - e^{ A { ( s - kh ) } } O_{k h} \bigr) \,ds
\Biggr\|_H^2 \Biggr]
\\[-0.2pt]
&&\qquad= \sum_{k=0}^{m-2} \mathbb{E} \biggl[ \biggl\| \int_{k h}^{( k + 1 ) h} e^{ A { ( m
h - s ) } } F'(X_{k h}) \bigl( O_s - e^{ A { ( s - kh ) } } O_{k h} \bigr) \,ds
\biggr\|_H^2 \biggr]
\\[-0.2pt]
&&\qquad\quad{}+ 2 \mathop{\sum_{k,\tilde{k} = 0}}_{k < \tilde{k}}^{m - 2} \mathbb{E}
\biggl[ \biggl\langle \int_{k h}^{( k + 1 ) h} e^{ A { ( m h - s ) } } F'(X_{k h})\\
&&\qquad\quad\hspace*{95.3pt}{}\times \bigl( O_s
- e^{ A { ( s - kh ) } } O_{k h} \bigr) \,ds,
\\[-0.2pt]
&&\qquad\quad\hspace*{61pt}\int_{\tilde{k} h}^{( \tilde{k} + 1 ) h} e^{ A { ( m h - s ) } }
F'(X_{\tilde{k} h}) \\[-0.2pt]
&&\qquad\quad\hspace*{95.4pt}{}\times\bigl( \mathbb{E} \bigl[ O_s - e^{ A { ( s - \tilde{k} h ) }
} O_{\tilde{k} h} | \mathcal{F}_{\tilde{k} h} \bigr] \bigr) \,ds
\biggr\rangle_H\biggr]
\end{eqnarray*}
and
%
%
\begin{eqnarray} \label{eqPart6}
&&\mathbb{E} \Biggl[ \Biggl\| \sum_{k = 0}^{m - 2} \int_{k h}^{( k + 1 ) h} e^{ A { (
m h - s ) } } F'(X_{k h}) \bigl( O_s - e^{ A { ( s - kh ) } } O_{k h} \bigr) \,ds
\Biggr\|_H^2 \Biggr]
\nonumber\\[-8pt]\\[-8pt]
&&\qquad= \sum_{k=0}^{m-2} \mathbb{E} \biggl[ \biggl\| \int_{k h}^{( k + 1 ) h} e^{ A { ( m
h - s ) } } F'(X_{k h}) \bigl( O_s - e^{ A { ( s - kh ) } } O_{k h} \bigr) \,ds
\biggr\|_H^2 \biggr]\nonumber
\end{eqnarray}
for every $ m \in\{ 0,1,\ldots, M \} $
and every $ M \in\mathbb{N} $ due to
Assumption~\ref{stochconv}.
Combining~(\ref{eqPart5}) and~(\ref{eqPart6}) then shows
\begin{eqnarray*}
&& \| X_{m h} - Y_m^M \|_{ L^2 ( \Omega; H ) }
\\[-2pt]
&&\qquad\leq
c \sum_{k = 0}^{m - 2} \int_{k h}^{( k + 1 ) h} \| ( X_s - O_s ) - (
X_{k h} - O_{k h} ) \|_{ L^2 ( \Omega; H ) } \,ds
\\[-2pt]
&&\qquad\quad{}
+ \Biggl( \sum_{k = 0}^{m - 2} \biggl\| \int_{k h}^{( k + 1 ) h} e^{ A { ( m h - s
) } } F'(X_{k h}) \bigl( O_s - e^{ A { ( s - kh ) } } O_{k h} \bigr) \,ds \biggr\|_{ L^2
( \Omega; H ) }^2 \Biggr)^{1/2}
\\[-2pt]
&&\qquad\quad{}
+ \sum_{k = 0}^{m - 2} \int_{k h}^{( k + 1 ) h} \bigl\| e^{ A { ( m h - s )
} } F'( X_{ k h } ) \bigl( \bigl( e^{ A { ( s - k h ) } } - I \bigr) O_{ k h } \bigr) \bigr\|_{
L^2 ( \Omega; H ) } \,ds
\\[-2pt]
&&\qquad\quad{}
+ c R \sum_{k = 0}^{m - 2} \int_{k h}^{( k + 1 ) h } \| O_s - O_{k h}
\|_{ L^4 ( \Omega; V ) }^2 \,ds + \sum_{k = 0}^{m - 2} \int_{k h}^{( k +
1 ) h } R \frac{ ( s - k h ) }{ ( m h - s ) } \,ds
\\[-2pt]
&&\qquad\quad{}
+ 2 R^2 M^{-1}
\end{eqnarray*}
for every $ m \in\{ 0,1,\ldots, M \} $
and every $ M \in\mathbb{N} $.
Hence, we obtain
\begin{eqnarray*}
&& \| X_{m h} - Y_m^M \|_{ L^2 ( \Omega; H ) }
\\[-2pt]
&&\qquad\leq
c R \sum_{k = 0}^{m - 2} \int_{k h}^{( k + 1 ) h} ( s - k h ) \,ds\\[-2pt]
&&\qquad\quad{} +
\sum_{k = 0}^{m - 2} \int_{k h}^{( k + 1 ) h } R \frac{ ( s - k h ) }{
( m h - (k+1)h ) } \,ds
\\[-2pt]
&&\qquad\quad{}
+ \Biggl( \sum_{k = 0}^{m - 2} \biggl\| \int_{k h}^{( k + 1 ) h} e^{ A { ( m h - s
) } } F'(X_{k h}) \bigl( O_s - e^{ A { ( s - kh ) } } O_{k h} \bigr) \,ds \biggr\|_{ L^2
( \Omega; H ) }^2 \Biggr)^{1/2}
\\[-2pt]
&&\qquad\quad{}
+ \sum_{k = 0}^{m - 2} \int_{k h}^{( k + 1 ) h} \bigl\| e^{ A { ( m h - s )
} } (-A)^{1/2} \bigr\|_{L(H)}
\\[-2pt]
&&\qquad\quad\hspace*{66.1pt}{}
\times\bigl\| (-A)^{- {1/2}} F'( X_{ k h } ) \bigl( \bigl( e^{ A { ( s - k h ) }
} - I \bigr) O_{ k h } \bigr) \bigr\|_{ L^2 ( \Omega; H ) } \,ds
\\[-2pt]
&&\qquad\quad{}
+ c R \sum_{k = 0}^{m - 2} \int_{k h}^{( k + 1 ) h } \| O_s - O_{k h}
\|_{ L^4 ( \Omega; V ) }^2 \,ds + 2 R^2 M^{-1}
\end{eqnarray*}
and
\begin{eqnarray*}
&& \| X_{m h} - Y_m^M \|_{ L^2 ( \Omega; H ) }
\\[-2pt]
&&\qquad\leq
\frac{1}{2} c R M h^2\\
&&\qquad\quad{} + R \sum_{k = 0}^{m - 2} \frac{ h }{ 2 ( m - k -1
) } + 2 R^2 M^{-1}
\\
&&\qquad\quad{}
+ \Biggl\{ \sum_{k = 0}^{m - 2} \biggl( \int_{k h}^{( k + 1 ) h} \bigl\| e^{ A { ( m h -
s ) } } F'(X_{k h}) \\
&&\qquad\quad\hspace*{82.2pt}{}\times\bigl( O_s - e^{ A { ( s - kh ) } } O_{k h} \bigr) \bigr\|_{ L^2 (
\Omega; H ) } \,ds \biggr)^{ 2 } \Biggr\}^{{1/2}}
\\
&&\qquad\quad{}
+ \sum_{k = 0}^{m - 2} \int_{k h}^{( k + 1 ) h} (mh-s)^{ - {1/2}}\\
&&\qquad\quad\hspace*{65.4pt}{}\times
\bigl\| (-A)^{- {1/2}} F'( X_{ k h } ) \bigl( \bigl( e^{ A { ( s - k h ) } } - I
\bigr) O_{ k h } \bigr) \bigr\|_{ L^2 ( \Omega; H ) } \,ds
\\
&&\qquad\quad{}
+ c R \sum_{k = 0}^{m - 2} \int_{k h}^{( k + 1 ) h } \bigl( R(s-kh)^\theta
\bigr)^2 \,ds
\end{eqnarray*}
for every $ m \in\{ 0,1,\ldots, M \} $
and every $ M \in\mathbb{N} $.
This yields
\begin{eqnarray*}
&& \| X_{m h} - Y_m^M \|_{ L^2 ( \Omega; H ) }
\\
&&\qquad\leq
\frac{1}{2} c R T h + \frac{1}{2} R h \Biggl( \sum_{k = 1}^{m - 1}
\frac{1}{k} \Biggr) + 2 R^2 M^{-1}\\
&&\qquad\quad{} + c R^3 \sum_{k = 0}^{m - 2} \int_{k h}^{(
k + 1 ) h } (s-kh)^{2 \theta} \,ds
\\
&&\qquad\quad{}
+ \Biggl\{ \sum_{k = 0}^{m - 2} \biggl( \int_{k h}^{( k + 1 ) h} \bigl\| F'(X_{k h}) \bigl(
O_s - e^{ A { ( s - kh ) } } O_{k h} \bigr) \bigr\|_{ L^2 ( \Omega; H ) } \,ds \biggr)^2
\Biggr\}^{{1/2}}
\\
&&\qquad\quad{}
+ \sum_{k = 0}^{m - 2} \int_{k h}^{( k + 1 ) h} \bigl(mh-(k+1)h\bigr)^{ -1/2 }
\\
&&\qquad\quad\hspace*{65.7pt}{}
\times\bigl\| (-A)^{- {1/2}} F'( X_{ k h } ) \bigl( \bigl( e^{ A { ( s - k h ) }
} - I \bigr) O_{ k h } \bigr) \bigr\|_{ L^2 ( \Omega; H ) } \,ds
\end{eqnarray*}
and
\begin{eqnarray*}
&& \| X_{m h} - Y_m^M \|_{ L^2 ( \Omega; H ) }
\\
&&\qquad\leq
\frac{1}{2} c R T^2 M^{-1} + \frac{1}{2} R^2 M^{-1} \Biggl( 1 + \sum_{k =
2}^{m - 1} \frac{1}{k} \Biggr)\\
&&\qquad\quad{} + 2 R^2 M^{-1} + c R^3 M h^{( 1 + 2 \theta)}
\\
&&\qquad\quad{}
+ \sqrt{h} \Biggl\{ \sum_{k = 0}^{m - 2} \int_{k h}^{( k + 1 ) h} \bigl\| F'(X_{k
h}) \bigl( O_s - e^{ A { ( s - kh ) } } O_{k h} \bigr) \bigr\|_{ L^2 ( \Omega; H ) }^2
\,ds \Biggr\}^{ {1/2} }
\\[-2pt]
&&\qquad\quad{}
+ \sqrt{T} \sum_{k = 0}^{m - 2} \frac{1}{(m-k-1)h}
\int_{k h}^{( k + 1 ) h} \bigl\| (-A)^{- {1/2}} F'( X_{ k h } )
\\[-2pt]
&&\qquad\quad\hspace*{48.4pt}\hspace*{104.5pt}{}
\times\bigl( \bigl( e^{ A { ( s - k h ) } } - I \bigr) O_{ k h } \bigr) \bigr\|_{ L^2 ( \Omega; H ) }
\,ds
\end{eqnarray*}
for every $ m \in\{ 0,1,\ldots, M \} $
and every $ M \in\mathbb{N} $.
Hence, we have
\begin{eqnarray*}
&& \| X_{m h} - Y_m^M \|_{ L^2 ( \Omega; H ) }
\\[-2pt]
&&\qquad\leq
\frac{1}{2} R^4 M^{-1} + \frac{1}{2} R^2 M^{-1} \biggl( 1 + \int_1^M
\frac{1}{s} \,ds \biggr) + 2 R^2 M^{-1} + c R^3 T h^{ 2 \theta}
\\[-2pt]
&&\qquad\quad{}
+ \sqrt{h} \Biggl\{ \sum_{k = 0}^{m - 2} \int_{k h}^{( k + 1 ) h} c^2 \bigl\| O_s
- e^{ A { ( s - kh ) } } O_{k h} \bigr\|_{ L^2 ( \Omega; H ) }^2 \,ds \Biggr\}^{
{1/2} }
\\[-2pt]
&&\qquad\quad{}
+ R \sum_{k = 0}^{m - 2} \frac{1}{(m-k-1)h}
\\[-2pt]
&&\qquad\quad\hspace*{39.1pt}{}
\times\int_{k h}^{( k + 1 ) h} \bigl\| (-A)^{- {1/2}} F'( X_{ k h } )
\bigl( \bigl( e^{ A { ( s - k h ) } } - I \bigr) O_{ k h } \bigr) \bigr\|_{ L^2 ( \Omega; H ) }
\,ds
\end{eqnarray*}
and
\begin{eqnarray*}
&& \| X_{m h} - Y_m^M \|_{ L^2 ( \Omega; H ) }
\\[-2pt]
&&\qquad\leq
\frac{1}{2} R^4 M^{-1} + \frac{1}{2} R^2 M^{-1} \bigl( 1 + \log(M) \bigr) + 2 R^2
M^{-1} + R^6 M^{ -2 \theta}
\\[-2pt]
&&\qquad\quad{}
+ \sqrt{T} c M^{- {1/2} } \Biggl\{ \sum_{k = 0}^{m - 2} \int_{k h}^{( k
+ 1 ) h} \bigl\| O_s - e^{ A { ( s - kh ) } } O_{k h} \bigr\|_{ L^2 ( \Omega; H )
}^2 \,ds \Biggr\}^{ {1/2} }
\\[-2pt]
&&\qquad\quad{}
+ \sum_{k = 0}^{m - 2} \frac{R}{(m-k-1)h}
\\[-2pt]
&&\qquad\quad\hspace*{29.2pt}{}
\times\int_{k h}^{( k + 1 ) h} \Bigl\| \sup_{ \| w \|_H \leq1 } \bigl| \bigl\langle w,
(-A)^{- {1/2}} F'( X_{ k h } ) \\[-2pt]
&&\qquad\quad\hspace*{121.5pt}{}\times\bigl( e^{ A ( s - k h ) } - I \bigr) O_{ k
h } \bigr\rangle_H \bigr| \Bigr\|_{ L^2 ( \Omega; \mathbb{R} ) } \,ds
\end{eqnarray*}
for every $ m \in\{ 0,1,\ldots, M \} $
and every $ M \in\mathbb{N} $.
This yields
\begin{eqnarray*}
&& \| X_{m h} - Y_m^M \|_{ L^2 ( \Omega; H ) } \\[-2pt]
&&\qquad\leq\biggl( \frac{1}{2} R^4 +
\frac{1}{2} R^2 + 2 R^2 + R^6 \biggr) \frac{( 1 + \log(M))}{ M^{ 2 \theta} }
\\
&&\qquad\quad{}
+ \sum_{k = 0}^{m - 2} \frac{R}{(m-k-1)h}
\\
&&\qquad\quad{}
\times\int_{k h}^{( k + 1 ) h} \Bigl\| \sup_{ \| w \|_H \leq1 } \bigl| \bigl\langle(
F'( X_{ k h } ) )^* (-A)^{- {1/2}} w,\\
&&\qquad\quad\hspace*{116pt}
\bigl( e^{ A ( s - k h ) } - I \bigr) O_{ k h } \bigr\rangle_H \bigr|
\Bigr\|_{ L^2 ( \Omega; \mathbb{R} ) } \,ds\\
&&\qquad\quad{}
+ R^2 M^{- {1/2} }
\Biggl\{ \sum_{k = 0}^{m - 2} \int_{k h}^{( k + 1 ) h} \bigl( \| O_s -
O_{kh} \|_{L^2 ( \Omega; H ) } \\
&&\qquad\quad\hspace*{78.4pt}\hspace*{42.5pt}{} + \bigl\| e^{A(s-kh)} O_{kh} - O_{kh} \bigr\|_{
L^2 ( \Omega; H ) } \bigr)^2 \,ds \Biggr\}^{ {1/2} }
\end{eqnarray*}
and
\begin{eqnarray*}
&& \| X_{m h} - Y_m^M \|_{ L^2 ( \Omega; H ) } \\
&&\qquad\leq4 R^6 \frac{( 1 +
\log(M))}{ M^{ 2 \theta} }
\\
&&\qquad\quad{}
+ \sum_{k = 0}^{m - 2} \frac{R}{(m-k-1)h} \\
&&\qquad\quad\hspace*{29pt}{}\times\int_{k h}^{( k + 1 ) h}
\Bigl\| \sup_{ \| w \|_H \leq1 } \bigl\| \bigl( F'( X_{ k h } ) \bigr)^* ( -A )^{-
{1/2}} w \bigr\|_{ D ( (-A)^{1/2} )} \\
&&\qquad\quad\hspace*{116.5pt}{}\times
\bigl\| \bigl( e^{ A ( s - k h ) } - I \bigr) O_{ k h }\bigr\|_{ D ( (-A)^{-{1/2}}
)} \Bigr\|_{ L^2 ( \Omega; \mathbb{R} ) } \,ds
\\
&&\qquad\quad{}
+ R^3 M^{- {1/2} }
\\
&&\qquad\quad\hspace*{10pt}{}
\times\Biggl\{ \sum_{k = 0}^{m - 2} \int_{k h}^{( k + 1 ) h} \bigl( \| O_s -
O_{kh} \|_{ L^2 ( \Omega; V ) }\\
&&\qquad\quad\hspace*{89pt}{} + \bigl\| \bigl( e^{A(s-kh)} - I \bigr) O_{kh} \bigr\|_{
L^2 ( \Omega; H ) } \bigr)^2 \,ds \Biggr\}^{1/2}
\end{eqnarray*}
for every $ m \in\{ 0,1,\ldots, M \} $
and every $ M \in\mathbb{N} $.
Using now condition~(\ref{assumpt2_2}) in
Assumption~\ref{drift} shows
%
%
\begin{eqnarray}\label{keyproof}
&& \| X_{m h} - Y_m^M \|_{ L^2 ( \Omega; H ) } \nonumber\\
&&\qquad\leq4 R^6 \frac{( 1 +
\log(M))}{ M^{ 2 \theta} } \nonumber
\\
&&\qquad\quad{}
+ \sum_{k = 0}^{m - 2} \frac{R}{(m-k-1)h}
\nonumber\\
&&\qquad\quad\hspace*{29.4pt}{}\times
\int_{k h}^{( k + 1 ) h}
\bigl\| c \bigl( 1 + \| X_{ k h } \|_{D ( (-A)^{1/2} )} \bigr)
\nonumber\\[-8pt]\\[-8pt]
&&\qquad\quad\hspace*{84pt}{}\times\bigl\| \bigl( e^{ A ( s - k h ) } - I \bigr) O_{ k h } \bigr\|_{D (
(-A)^{-{1/2}} )} \bigr\|_{ L^2 ( \Omega; \mathbb{R} ) } \,ds
\nonumber\\
&&\qquad\quad{}
+ R^3 M^{- {1/2} }
\nonumber
\\
&&\qquad\quad\hspace*{10pt}{}
\times\Biggl\{ \sum_{k = 0}^{m - 2} \int_{k h}^{( k + 1 ) h} \bigl( R (s -
kh)^\theta\nonumber\\
&&\qquad\quad\hspace*{87pt}{}+ (s-kh)^\gamma\| O_{kh} \|_{ L^2 ( \Omega; { D (
(-A)^\gamma) } ) } \bigr)^2 \,ds \Biggr\}^{ {1/2} }\nonumber
\end{eqnarray}
and therefore
\begin{eqnarray*}
&& \| X_{m h} - Y_m^M \|_{ L^2 ( \Omega; H ) } \\
&&\qquad\leq4 R^6 \frac{( 1 +
\log(M))}{ M^{ 2 \theta} }
\\
&&\qquad\quad{}
+ \sum_{k = 0}^{m - 2} \frac{c R}{(m-k-1)h}
%
\int_{k h}^{( k + 1 ) h} \bigl\| 1 + \| X_{ k h } \|_{ D ( (-A)^{1/2}
) } \bigr\|_{ L^4( \Omega; \mathbb{R} ) }
\\
&&\qquad\quad\hspace*{130.2pt}{}
\times\bigl\| \bigl( e^{ A ( s - k h ) } - I \bigr) O_{ k h } \bigr\|_{ L^4( \Omega; D (
(-A)^{-{1/2}} ) ) } \,ds
\\
&&\qquad\quad{}
+ R^3 M^{- {1/2} } \Biggl\{ \sum_{k = 0}^{m - 2} \int_{k h}^{( k + 1 )
h} \bigl( R h^\theta+ R h^\theta T^{(\gamma- \theta)} \bigr)^2 \,ds \Biggr\}^{
{1/2} }
\end{eqnarray*}
for every $ m \in\{ 0,1,\ldots, M \} $
and every $ M \in\mathbb{N} $.
Hence, we obtain
\begin{eqnarray*}
&& \| X_{m h} - Y_m^M \|_{ L^2 ( \Omega; H ) }
\\
&&\qquad\leq
4 R^6 \frac{( 1 + \log(M))}{ M^{ 2 \theta} }\\
&&\qquad\quad{} + \sum_{k = 0}^{m - 2}
\frac{c R}{(m-k-1)h} \bigl( 1 + \| X_{ k h } \|_{L^4( \Omega; D(
(-A)^{1/2} ) ) } \bigr)
\\
&&\qquad\quad\hspace*{29.7pt}{}
\times\int_{k h}^{( k + 1 ) h} \bigl\| \bigl( e^{ A ( s - k h ) } - I \bigr) O_{ k h }
\bigr\|_{ L^4( \Omega; D ( (-A)^{-{1/2}} ) ) } \,ds
\\
&&\qquad\quad{}
+ R^3 M^{- {1/2} } \Biggl\{ \sum_{k = 0}^{m - 2} \int_{k h}^{( k + 1 )
h} ( 2 R^2 h^\theta)^2 \,ds \Biggr\}^{ {1/2} }
\end{eqnarray*}
and
\begin{eqnarray*}
&&\| X_{m h} - Y_m^M \|_{ L^2 ( \Omega; H ) }
\\
&&\qquad\leq4 R^6 \frac{( 1 + \log(M))}{ M^{ 2 \theta} }
\\
&&\qquad\quad{}
+ 2 R^5 M^{- {1/2} } \Biggl( \sum_{k = 0}^{m - 2} \int_{k h}^{( k + 1 )
h} h^{ 2 \theta} \,ds \Biggr)^{ {1/2} }\\
&&\qquad\quad{} + \sum_{k = 0}^{m - 2} \frac{ c R
(1 + R) }{(m-k-1)h}
\\
&&\qquad\quad\hspace*{29.2pt}{}
\times\int_{k h}^{( k + 1 ) h} \bigl\| (-A)^{ - ( \gamma+ {1/2} ) } \bigl(
e^{ A ( s - k h ) } - I \bigr) \bigr\|_{ L ( H ) }\\
&&\qquad\quad\hspace*{77pt}{}\times \| O_{ k h } \|_{ L^4( \Omega;
D( (-A)^\gamma) ) } \,ds
\end{eqnarray*}
for every $ m \in\{ 0,1,\ldots, M \} $
and every $ M \in\mathbb{N} $.
This yields
\begin{eqnarray*}
&&\| X_{m h} - Y_m^M \|_{ L^2 ( \Omega; H ) } \\
&&\qquad\leq4 R^6 \frac{( 1 +
\log(M))}{ M^{ 2 \theta} } + 2 R^5 M^{- {1/2} } \bigl( M h^{
(1+2\theta) } \bigr)^{1/2}
\\
&&\qquad\quad{}+ \sum_{k = 0}^{m - 2} \frac{ 2 c R^3 }{(m-k-1)h} \int_{k h}^{( k + 1 )
h} \bigl\| (-A)^{ - ( \gamma+ {1/2} ) } \bigl( e^{ A ( s - k h ) } - I \bigr)
\bigr\|_{ L ( H ) } \,ds
\end{eqnarray*}
and hence
\begin{eqnarray*}
&&\| X_{m h} - Y_m^M \|_{ L^2 ( \Omega; H ) } \\
&&\qquad\leq4 R^6 \frac{( 1 +
\log(M))}{ M^{ 2 \theta} } + 2 R^5 \sqrt{T} M^{- {1/2} } h^\theta
\\
&&\qquad\quad{} + \sum_{k = 0}^{m - 2} \frac{ 2 R^4 }{(m-k-1)h} \int_{k h}^{( k + 1 )
h} \bigl\| (-A)^{ ( {1/2} - \gamma) } \bigr\|_{ L ( H ) } \\
&&\qquad\quad\hspace*{130.2pt}{}\times\bigl\| A^{-1} \bigl( e^{ A
( s - k h ) } - I \bigr) \bigr\|_{ L ( H ) } \,ds
\end{eqnarray*}
for every $ m \in\{ 0,1,\ldots, M \} $
and every $ M \in\mathbb{N} $.
Therefore, we have
\begin{eqnarray*}
&&\| X_{m h} - Y_m^M \|_{ L^2 ( \Omega; H ) } \\
&&\qquad\leq4 R^6 \frac{( 1 +
\log(M))}{ M^{ 2 \theta} } + 2 R^6 M^{- ( {1/2} + \theta) }
\\
&&\qquad\quad{}+ \sum_{k = 0}^{m - 2} \frac{ 2 R^4 }{(m-k-1)h} \int_{k h}^{( k + 1 )
h} \biggl( \frac{1}{\lambda_1} \biggr)^{ ( \gamma- {1/2} ) } ( s - k h ) \,ds
\end{eqnarray*}
and, finally,
%
%
\begin{eqnarray} \label{nrA}
\| X_{m h} - Y_m^M \|_{ L^2 ( \Omega; H ) } &\leq&6 R^6 \frac{( 1 +
\log(M))}{ M^{ 2 \theta} } + \sum_{k = 0}^{m - 2} \frac{ R^5 h
}{(m-k-1)} \nonumber
\\[-2pt]
&\leq&
6 R^6 \frac{( 1 + \log(M))}{ M^{ 2 \theta} } + R^6 M^{-1} \Biggl( \sum_{k =
1}^{M} \frac{1}{k} \Biggr) \nonumber\\[-9pt]\\[-9pt]
&\leq&
6 R^6 \frac{( 1 + \log(M))}{ M^{ 2 \theta} } + R^6 M^{-2 \theta} \biggl( 1 +
\int_1^{M} \frac{1}{s} \,ds \biggr) \nonumber\\[-2pt]
&=& 7 R^6 \frac{( 1 + \log(M))}{ M^{ 2
\theta} }\nonumber
\end{eqnarray}
for every $ m \in\{ 0,1,\ldots, M \} $
and every $ M \in\mathbb{N} $.
%
\subsubsection{Spatial discretization error}\label{spat_error}
Due to~(\ref{NumScheme1}),
we obtain
\begin{eqnarray*}
&& \| Y_m^M - P_N( Y_m^M ) \|_{ L^2 ( \Omega; H ) }
\\[-2pt]
&&\qquad= \Biggl\| e^{ A { m h } } \bigl( \xi- P_N( \xi) \bigr)\\[-2pt]
&&\qquad\quad\hspace*{3.1pt}{} + h \Biggl( \sum_{ k = 0 }^{ m - 1 }
\bigl( e^{ A { ( m h - k h ) } } - P_N e^{ A { ( m h - k h ) } } \bigr) F( X_{ k
h } ) \Biggr)
\\[-2pt]
&&\qquad\quad\hspace*{120.6pt}{}
+ O_{ m h } - P_N( O_{ m h } ) \Biggr\|_{ L^2 ( \Omega; H ) }
\\[-2pt]
&&\qquad\leq
\bigl\| e^{ A { m h } } \bigl( \xi- P_N( \xi) \bigr) \bigr\|_{ L^2 ( \Omega; H ) } + \| O_{
m h } - P_N( O_{ m h } ) \|_{ L^2 ( \Omega; H ) }
\\[-2pt]
&&\qquad\quad{} + \Biggl\| h \Biggl( \sum_{ k = 0 }^{ m - 1 } \bigl( e^{ A { ( m h - k h ) } } - P_N
e^{ A { ( m h - k h ) } } \bigr) F( X_{ k h } ) \Biggr) \Biggr\|_{ L^2 ( \Omega; H ) }
\end{eqnarray*}
and
\begin{eqnarray*}
&& \| Y_m^M - P_N( Y_m^M ) \|_{ L^2 ( \Omega; H ) }
\\[-2pt]
&&\qquad\leq
\| \xi- P_N( \xi) \|_{ L^2 ( \Omega; H ) } + \| O_{ m h } - P_N( O_{ m
h } ) \|_{ L^2 ( \Omega; H ) }
\\[-2pt]
&&\qquad\quad{} + h \Biggl( \sum_{ k = 0 }^{ m - 1 } \bigl\| e^{ A { ( m h - k h ) } } - P_N e^{
A { ( m h - k h ) } } \bigr\|_{ L ( H ) } \| F( X_{ k h } ) \|_{ L^2 (
\Omega; H ) } \Biggr)
\\[-2pt]
&&\qquad\leq
\| (I-P_N) \xi\|_{ L^2 ( \Omega; H ) } + \| (I-P_N) O_{mh} \|_{ L^2 (
\Omega; H ) }
\\[-2pt]
&&\qquad\quad{}+ R h \Biggl( \sum_{ k = 0 }^{ m - 1 } \bigl\| ( I - P_N ) e^{ A { ( m h - k h ) }
} \bigr\|_{ L(H) } \Biggr)
\end{eqnarray*}
and hence
\begin{eqnarray*}
&&\| Y_m^M - P_N( Y_m^M ) \|_{ L^2 ( \Omega; H ) }\\[-2pt]
&&\qquad\leq\| (-A)^{ - \gamma} ( I - P_N ) \|_{ L ( H ) } \| (-A)^\gamma\xi
\|_{ L^2 ( \Omega; H ) }
\\[-2pt]
&&\qquad\quad{} + \| (-A)^{ - \gamma} ( I - P_N ) \|_{ L ( H ) } \| (-A)^\gamma
O_{mh} \|_{ L^2 ( \Omega; H ) }
\\[-2pt]
&&\qquad\quad{}+ R h \Biggl( \sum_{ k = 0 }^{ m - 1 } \| (-A)^{ - \gamma} ( I - P_N ) \|_{
L ( H ) } \bigl\| (-A)^\gamma e^{ A { ( m h - k h ) } } \bigr\|_{ L ( H ) } \Biggr)
\\[-2pt]
&&\qquad\leq(\lambda_N)^{-\gamma} \bigl( \| (-A)^\gamma\xi\|_{ L^2 ( \Omega; H )
} + \| (-A)^\gamma O_{mh} \|_{ L^2 ( \Omega; H ) } \bigr)
\\[-2pt]
&&\qquad\quad{} + R h \Biggl( \sum_{ k = 0 }^{ m - 1 } (\lambda_N)^{-\gamma} \bigl\| (-A)^\gamma
e^{ A { ( m h - k h ) } } \bigr\|_{ L ( H ) } \Biggr)
\end{eqnarray*}
for every $ m \in\{ 0,1,\ldots, M \} $
and every $ M \in\mathbb{N} $.
Therefore, we have
\begin{eqnarray*}
&& \| Y_m^M - P_N( Y_m^M ) \|_{ L^2 ( \Omega; H ) } \nonumber
\\[-2pt]
&&\qquad
\leq( \lambda_N )^{-\gamma} \bigl( \| \xi\|_{ L^2 ( \Omega; { D (
(-A)^\gamma) } ) } + \| O_{mh} \|_{ L^2( \Omega; { D ( (-A)^\gamma) } )
} \bigr) \nonumber
\\[-2pt]
&&\qquad\quad{}
+ R h (\lambda_N)^{ - \gamma} \Biggl( \sum_{ k = 0 }^{ m - 1 } \frac{1} { (m
h - k h)^\gamma} \bigl\| \bigl(-A (m h - k h) \bigr)^\gamma e^{ A { ( m h - k h ) } }
\bigr\|_{ L ( H ) } \Biggr) \nonumber
\\[-2pt]
&&\qquad\leq
2 R ( \lambda_N )^{-\gamma} + R h^{(1-\gamma)} (\lambda_N)^{-\gamma} \Biggl(
\sum_{ k = 0 }^{ m - 1 } \frac{1} { (m - k)^\gamma} \Bigl( \sup_{ x>0 }
x^\gamma e^{-x} \Bigr) \Biggr) \nonumber
\\[-2pt]
&&\qquad\leq
2 R ( \lambda_N )^{-\gamma} + R h^{(1-\gamma)} (\lambda_N)^{-\gamma} \Biggl(
\sum_{ k = 1 }^{ m } \frac{1}{k^\gamma} \Biggr)
\\[-2pt]
&&\qquad\leq
2 R ( \lambda_N )^{-\gamma} + R h^{(1-\gamma)} (\lambda_N)^{-\gamma} \Biggl(
1 + \sum_{ k = 2 }^{ m } \frac{1}{k^\gamma} \Biggr) \nonumber\vspace*{-2pt}
\end{eqnarray*}
and
%
%
\begin{eqnarray} \label{nrB}
&&
\| Y_m^M - P_N( Y_m^M ) \|_{ L^2 ( \Omega; H ) } \nonumber\\[-3pt]
&&\qquad\leq R (
\lambda_N )^{-\gamma} \biggl( 2 + h^{(1-\gamma)} \biggl( 1 + \int_1^M
\frac{1}{s^\gamma} \,ds \biggr) \biggr) \nonumber
\\[-3pt]
&&\qquad= R ( \lambda_N )^{-\gamma} \biggl( 2 + h^{(1-\gamma)} \biggl( 1 + \biggl[
\frac{s^{(1-\gamma)}} {(1-\gamma)} \biggr]_{s=1}^{s=M} \biggr) \biggr)
\\[-3pt]
&&\qquad= R ( \lambda_N )^{-\gamma} \biggl( 2 + h^{(1-\gamma)} \biggl( 1 +
\frac{M^{(1-\gamma)}} {(1-\gamma)} - \frac{1} {(1-\gamma)} \biggr) \biggr)
\nonumber\\[-3pt]
&&\qquad\leq R ( \lambda_N )^{-\gamma} \biggl( 2 + \frac{T^{(1-\gamma)}}
{(1-\gamma)} \biggr) \leq3 R^3 ( \lambda_N )^{-\gamma} \nonumber
\end{eqnarray}
for every $ m \in\{ 0,1, \ldots, M
\} $ and every $ M \in\mathbb{N} $.\vadjust{\goodbreak}
%
\subsubsection{Lipschitz estimates} \label{refLipschitzEstimate}
Note that $ Y_m^{N, M} \dvtx\Omega\rightarrow V $
satisfies
%
%
\begin{equation}\label{lip_1}\qquad
Y_m^{N,M} = e^{ A { m h } } ( P_N( \xi) ) + h \Biggl( \sum_{ k = 0 }^{ m - 1
} P_N e^{ A { ( m h - k h ) } } F( Y_k^{N,M} ) \Biggr) + P_N( O_{ m h } )
\end{equation}
for every $ m \in\{ 0, 1, \ldots, M \} $
and every $ N, M \in\mathbb{N} $.
Indeed, in the case $ m = 0 $ we have
\begin{eqnarray*}
Y_0^{N,M} &=& P_N( \xi) + P_N( O_0 )
\\
& =& e^{ A 0 } ( P_N( \xi) ) + h \Biggl( \sum_{ k = 0 }^{ - 1 } P_N e^{ A { (
0 - k h ) } } F( Y_k^{N,M} ) \Biggr) + P_N( O_0 )
\end{eqnarray*}
for every $ N, M \in\mathbb{N} $.
Moreover, if~(\ref{lip_1}) holds for one
$ m \in\{0, 1, \ldots, M - 1 \} $,
then we obtain
\begin{eqnarray*}
Y_{m+1}^{N, M} &=& e^{ A h } \bigl( Y_m^{N,M} + h \cdot( P_N F ) ( Y_m^{N,M}
) \bigr) + P_N\bigl( O_{(m+1)h} - e^{ A h } O_{ m h } \bigr)
\\
&=& e^{ A h } Y_m^{N,M} + h \cdot P_N e^{A h} F( Y_m^{N,M} ) + P_N\bigl(
O_{(m+1)h} \bigr) - e^{ A h } P_N( O_{ m h } )
\\
&=& e^{A h}\bigl( Y_m^{N,M} - P_N( O_{ m h } ) \bigr) + h \cdot P_N e^{A h} F(
Y_m^{N,M} ) + P_N\bigl( O_{(m+1)h} \bigr)
\end{eqnarray*}
and
\begin{eqnarray*}
Y_{m+1}^{N, M} &=& e^{A h}\Biggl( e^{ A { m h } } ( P_N( \xi) ) + h \Biggl( \sum_{ k
= 0 }^{ m - 1} P_N e^{ A {( m h - k h )} } F( Y_k^{N,M} ) \Biggr) \Biggr)
\\
&&{} + h \cdot P_N e^{A h} F( Y_m^{N,M} ) + P_N\bigl( O_{(m+1)h} \bigr)
\\
&=& e^{ A { ( m + 1 ) h } } ( P_N( \xi) ) + h \Biggl( \sum_{ k = 0 }^{ m - 1}
P_N e^{ A {( ( m + 1 ) h - k h )} } F( Y_k^{N,M} ) \Biggr)
\\
&&{} + h \cdot P_N e^{A h} F( Y_m^{N,M} ) + P_N\bigl( O_{(m+1)h} \bigr)
\\
&=& e^{ A { ( m + 1 ) h } } ( P_N( \xi) ) + h \Biggl( \sum_{ k = 0 }^{ m } P_N
e^{ A {( ( m + 1 ) h - k h )} } F( Y_k^{N,M} ) \Biggr)
\\
&&{}
+ P_N\bigl( O_{(m+1)h} \bigr)
\end{eqnarray*}
for every $ N, M \in\mathbb{N} $,
which shows~(\ref{lip_1}) by induction.
In the next step,~(\ref{lip_1}) yields
\begin{eqnarray*}
&& P_N( Y_m^{N} ) - Y_m^{N, M}
\\
&&\qquad
= h \Biggl( \sum_{ k = 0 }^{ m - 1 } P_N e^{ A { ( m h - k h ) } } F( X_{ k h
} ) \Biggr)
\\
&&\qquad\quad{} - h \Biggl( \sum_{ k = 0 }^{ m - 1 } P_N e^{ A { ( m h - k h ) } } F( Y_{ k
}^{N, M} ) \Biggr)
\\
&&\qquad = h \Biggl( \sum_{ k = 0 }^{ m - 1 } P_N e^{ A { ( m h - k h ) } } \bigl( F( X_{
k h } ) - F( Y_{ k }^{N, M} ) \bigr) \Biggr)
\end{eqnarray*}
for every $ m \in\{ 0, 1, \ldots, M \} $
and every $ N, M \in\mathbb{N} $.
Therefore, we obtain
%
%
\begin{eqnarray} \label{nrC}
&& \| P_N( Y_m^{N} ) - Y_m^{N, M} \|_{ L^2 ( \Omega; H ) } \nonumber
\\
&&\qquad\leq
h \sum_{ k = 0 }^{ m - 1 } \bigl\| P_N e^{ A { ( m h - k h ) } } \bigl( F( X_{ k
h } ) - F( Y_{ k }^{N, M} ) \bigr) \bigr\|_{ L^2 ( \Omega; H ) } \nonumber
\\
&&\qquad
\leq h \sum_{ k = 0 }^{ m - 1 } \bigl( \bigl\| P_N e^{ A { ( m h - k h ) } } \bigr\|_{
L( H ) } \| F( X_{ k h } ) - F( Y_{ k }^{N, M} ) \|_{ L^2 ( \Omega; H )
} \bigr)
\\
&&\qquad
\leq h \sum_{ k = 0 }^{ m - 1 } \| F( X_{ k h } ) - F( Y_{ k }^{N, M} )
\|_{ L^2 ( \Omega; H ) } \nonumber\\
&&\qquad\leq c h \sum_{k = 0}^{m - 1} \| X_{kh} -
Y_k^{N, M} \|_{ L^2 ( \Omega; H ) }\nonumber
\end{eqnarray}
for every $ m \in\{ 0, 1, \ldots, M \} $
and every $ N, M \in\mathbb{N} $.
Combining~(\ref{nrA}),~(\ref{nrB}) and~(\ref{nrC}) finally yields
\begin{eqnarray*}
&&\| X_{m h} - Y_m^{ N, M } \|_{ L^2 ( \Omega; H ) } \\
&&\qquad\leq\| X_{m h} -
Y_m^M \|_{ L^2 ( \Omega; H ) }
\\
&&\qquad\quad{}
+ \| Y_m^M - P_N ( Y_m^M ) \|_{ L^2 ( \Omega; H ) } + \| P_N ( Y_m^M )
- Y_m^{ N, M } \|_{ L^2 ( \Omega; H ) }
\\
&&\qquad\leq
7 R^6 \frac{ ( 1 + \log(M) ) }{M^{2 \theta}} + 3 R^3 \frac{1}{
(\lambda_N)^{\gamma} } + c h \sum_{k = 0}^{m - 1} \| X_{kh} - Y_k^{N,
M} \|_{ L^2 ( \Omega; H ) }
\end{eqnarray*}
for every $ m \in\{ 0, 1, \ldots, M \} $
and every $ N, M \in\mathbb{N} $.
Hence, Gronwall's lemma yields
%
%
\begin{eqnarray}
\| X_{m h} - Y_m^{ N, M } \|_{ L^2 ( \Omega; H ) } &\leq&\biggl( 7 R^6
\frac{ ( 1 + \log(M) ) }{M^{2 \theta}} + 3 R^3 \frac{1}{( \lambda_N
)^{\gamma}} \biggr) e^{c T} \nonumber
\\
&\leq& \biggl( 7 R^6 \frac{ ( 1 + \log(M) ) }{M^{2 \theta}} + 7 R^6 \frac{1}{(
\lambda_N )^{\gamma}} \biggr) e^{ c T }
\\
&=& ( e^{ c T } 7 R^6 ) \biggl( \frac{ ( 1 + \log(M) ) }{M^{2 \theta}} +
\frac{1}{ ( \lambda_N )^{\gamma} } \biggr)\nonumber
\end{eqnarray}
for every $ m \in\{ 0, 1, \ldots, M \} $
and every $ N, M \in\mathbb{N} $,
which shows the assertion.
%
\subsection{\texorpdfstring{Properties of the SPDE (\protect\ref{existence})}
{Properties of the SPDE (1)}}\label{secprop}

\mbox{}

\begin{pf*}{Proof of Lemma~\ref{existence}}
A standard application of Banach's fix point theorem (see, e.g.,
Section 7.1 in~\cite{dz92}) yields the existence of a unique
adapted stochastic process $X\dvtx[0,T] \times\Omega\rightarrow V $
with continuous sample paths which fulfills~(\ref{eqSolution}).
%
Moreover, we have
%
%
\begin{equation} \label{proof1}
\int_0^t e^{A(t-s)} F(X_s(\omega)) \,ds \in{ D ( (-A)^\gamma) }
\end{equation}
for every $t \in[0,T] $ and every $\omega\in\Omega$, since
\begin{eqnarray*}
&& \int_0^t \bigl\| (-A)^\gamma e^{A(t-s)} F(X_s(\omega)) \bigr\|_H \,ds
\\[-2pt]
&&\qquad\leq
\int_0^t \bigl\| (-A)^\gamma e^{A(t-s)} \bigr\|_{ L ( H ) } \| F(X_s(\omega))
\|_H \,ds \\[-2pt]
&&\qquad\leq\int_0^t (t-s)^{-\gamma} \| F(X_s(\omega)) \|_H \,ds
\\[-2pt]
&&\qquad\leq
\int_0^t (t - s)^{-\gamma} \bigl( c \| X_s(\omega) \|_H + \| F(0) \|_H \bigr) \,ds
\\[-2pt]
&&\qquad\leq
\biggl( \int_0^t s^{-\gamma} \,ds \biggr) \Bigl( c \Bigl( \sup_{0 \leq s \leq T} \| X_s(\omega)
\|_H \Bigr) + \| F(0) \|_H \Bigr)
\end{eqnarray*}
and
\begin{eqnarray*}
&& \int_0^t \bigl\| (-A)^\gamma e^{A(t-s)} F(X_s(\omega)) \bigr\|_H \,ds
\\[-2pt]
&&\qquad\leq
\biggl[ \frac{s^{(1-\gamma)}} {(1-\gamma)} \biggr]_{s = 0}^{s = T} \Bigl( c \Bigl( {\sup_{0
\leq s \leq T}} \| X_s(\omega) \|_H \Bigr) + \| F(0) \|_H \Bigr)
\\[-2pt]
&&\qquad\leq
\frac{T^{(1-\gamma)}} {(1-\gamma)} \Bigl( c \Bigl( {\sup_{0 \leq s \leq T}} \|
X_s(\omega) \|_H \Bigr) + \| F(0) \|_H \Bigr) < \infty
\end{eqnarray*}
holds for every $ t \in[0,T] $
and every $\omega\in\Omega$.
Assumptions~\ref{stochconv},~\ref{initial} and
(\ref{proof1}) hence imply $X_t(\omega) \in
{ D ( (-A)^\gamma) }$ for every
$ t \in[0,T] $ and every $\omega\in\Omega$.
Furthermore, we have
\begin{eqnarray*}
&& \| (-A)^\gamma X_t \|_H
\\[-2pt]
&&\qquad\leq
\| (-A)^\gamma e^{At} \xi\|_H + \int_0^t \bigl\| (-A)^\gamma e^{A(t-s)}
F(X_s) \bigr\|_H \,ds + \| (-A)^\gamma O_t \|_H
\\[-2pt]
&&\qquad\leq
\| (-A)^\gamma\xi\|_H + \int_0^t \bigl\| (-A)^\gamma e^{A(t-s)} \bigr\|_{ L ( H
) } \| F(X_s) \|_H \,ds \\[-2pt]
&&\qquad\quad{} + {\sup_{0 \leq s \leq T}} \| (-A)^\gamma O_s \|_H
\\[-2pt]
&&\qquad\leq
\| (-A)^\gamma\xi\|_H + \int_0^t (t-s)^{-\gamma} \| F(X_s) \|_H \,ds +
{\sup_{0 \leq s \leq T}} \| (-A)^\gamma O_s \|_H
\\[-2pt]
&&\qquad\leq
\Bigl( \| (-A)^\gamma\xi\|_H + {\sup_{0 \leq s \leq T}} \| (-A)^\gamma O_s
\|_H \Bigr)
\\[-2pt]
&&\qquad\quad{}
+ \int_0^t (t-s)^{-\gamma} \bigl( c \| X_s \|_H + \| F(0) \|_H \bigr) \,ds
\end{eqnarray*}
for every $ t \in[0,T] $.
This yields
\begin{eqnarray*}
&&\| (-A)^\gamma X_t \|_H \\[-2pt]
&&\qquad\leq c \int_0^t (t-s)^{-\gamma} \| X_s \|_H
\,ds
\\[-2pt]
&&\qquad\quad{}
+ \biggl( \| (-A)^\gamma\xi\|_H + {\sup_{0 \leq s \leq T}} \| (-A)^\gamma O_s
\|_H + \| F(0) \|_H \biggl( \int_0^t s^{-\gamma} \,ds \biggr) \biggr)
\\[-2pt]
&&\qquad\leq
\biggl( \| (-A)^\gamma\xi\|_H + {\sup_{0 \leq s \leq T}} \| (-A)^\gamma O_s
\|_H + \frac{ T^{ (1 - \gamma)} \| F(0) \|_H } {(1-\gamma)} \biggr)
\\[-2pt]
&&\qquad\quad{}
+ c \| (-A)^{-\gamma} \|_{ L ( H ) } \int_0^t (t-s)^{-\gamma} \|
(-A)^\gamma X_s \|_H \,ds
\end{eqnarray*}
for every $ t \in[0,T] $.
Hence, Lemma 7.1.1 in~\cite{h81} shows
\begin{eqnarray*}
&&{\sup_{0 \leq t \leq T} }\| (-A)^\gamma X_t \|_H \\[-2pt]
&&\qquad\leq E_{(1 - \gamma)} \bigl(
T \bigl( c \| (-A)^{-\gamma} \|_{ L ( H ) } \Gamma( 1 - \gamma) \bigr)^{
{1}/{(1 - \gamma)} } \bigr)
\\[-2pt]
&&\qquad\quad\hspace*{0pt}{}\times\biggl( \| (-A)^\gamma\xi\|_H + {\sup_{0 \leq s \leq T}} \| (-A)^\gamma
O_s \|_H + \frac{ T^{(1 - \gamma)} \| F(0) \|_H } {(1-\gamma)} \biggr)
\end{eqnarray*}
and therefore
\begin{eqnarray*}
&& \Bigl\| \sup_{0 \leq t \leq T} \| (-A)^\gamma X_t \|_H \Bigr\|_{ L^4( \Omega;
\mathbb{R} ) } \\[-2pt]
&&\qquad\leq E_{(1 - \gamma)} \bigl( T \bigl( c \| (-A)^{-\gamma} \|_{ L (
H ) } \Gamma( 1 - \gamma) \bigr)^{ {1}/{(1 - \gamma)} } \bigr)
\\[-2pt]
&&\qquad\quad\hspace*{0pt}{}
\times\biggl( \| (-A)^\gamma\xi\|_{ L^4 ( \Omega; H ) } + \Bigl\| \sup_{0 \leq s
\leq T} \| (-A)^\gamma O_s \|_H \Bigr\|_{ L^4 ( \Omega; \mathbb{R} ) }\\[-2pt]
&&\qquad\quad\hspace*{160pt}{} +
\frac{ T^{(1 - \gamma)} \| F(0) \|_H } {(1-\gamma)} \biggr) < \infty,
\end{eqnarray*}
which shows the assertion.
Here $ E_{(1-\gamma)} \dvtx[0,\infty)
\rightarrow[0,\infty) $ is
given by
\[
E_{(1-\gamma)}(x) := \sum_{ n = 0 }^{ \infty} \frac{ x^{ (n ( 1 -
\gamma)) } }{ \Gamma( n ( 1 - \gamma) + 1 ) }
\]
for every $ x \in[0,\infty) $ where
$ \Gamma\dvtx(0,\infty) \rightarrow(0,\infty) $
is the Gamma function.\vadjust{\goodbreak}
\end{pf*}
\begin{lemma} \label{lemmaRef4}
Let Assumptions~\ref{semigroup}--\ref{initial} be fulfilled.
Then we have
\[
\| e^{At_2} - e^{At_1} \|_{ L ( H ) } \leq\frac{(t_2 - t_1)}{t_1}
\]
for every $ t_1, t_2 \in(0,T] $ with $t_1 \leq t_2$.
\end{lemma}
\begin{pf}
By definition, we have
\begin{eqnarray*}
&&
\| e^{At_2} - e^{At_1} \|_{ L ( H ) } \\
&&\qquad= \bigl\| \bigl( e^{A (t_2 - t_1) } - I \bigr)
e^{At_1} \bigr\|_{ L ( H ) }\\
&&\qquad\leq \bigl\| A^{-1} \bigl( e^{A (t_2 - t_1) } - I \bigr) \bigr\|_{ L ( H ) } \| A e^{At_1}
\|_{ L ( H ) }
\\
&&\qquad= \bigl\| \bigl( A (t_2 - t_1) \bigr)^{-1} \bigl( e^{A(t_2 - t_1)} - I \bigr) \bigr\|_{L(H)} \\
&&\qquad\quad\hspace*{0pt}{}\times\| A t_1
e^{A t_1} \|_{L(H)} \frac{(t_2 - t_1)}{t_1}
\\
&&\qquad\leq\biggl( \sup_{x \in(0, \infty)} \frac{( 1 - e^{-x} )} {x} \biggr) \Bigl( \sup_{x
\in(0, \infty)} x e^{-x} \Bigr) \frac{(t_2 - t_1)}{t_1} \\
&&\qquad\leq\frac{(t_2 -
t_1)}{t_1}
\end{eqnarray*}
for every $ t_1, t_2 \in(0,T] $ with $t_1 < t_2$.
\end{pf}
\begin{lemma} \label{lemmaRef2}
Let Assumptions~\ref{semigroup}--\ref{initial} be fulfilled.
Then we obtain
\[
\sup_{0 \leq t_1 < t_2 \leq T} \frac{ \| X_{t_2} - X_{t_1} \|_{ L^2 (
\Omega; H ) } } { ( t_2 - t_1 )^\theta} < \infty,
\]
where $\theta\in(0, \frac{1}{2} ]$
is given in Assumption~\ref{stochconv} and where
$X \dvtx\Omega\times[0,T] \rightarrow D
( (-A)^\gamma) $ is the solution of
the SPDE~(\ref{eqSolution}).
\end{lemma}
\begin{pf}
First, let $R \in[0, \infty) $ be the real number given by
\begin{eqnarray*}
R &:=& \| \xi\|_{ L^2 ( \Omega; D(A) ) } + \sup_{ t \in[0,T] } \| F(X_t)
\|_{ L^2 ( \Omega; H ) } \\
&&{}+ {\sup_{0 \leq t_1 < t_2 \leq T}} \biggl( \frac{ \|
O_{t_2} - O_{t_1} \|_{ L^2 ( \Omega; H ) } } { ( t_2 - t_1 )^\theta}
\biggr),
\end{eqnarray*}
which is finite due to Assumptions~\ref{semigroup}--\ref{initial}.
Then we have
%
%
\begin{eqnarray} \label{eqPart1}
&& \| e^{ A {t_2} } \xi- e^{ A {t_1} } \xi\|_{ L^2 ( \Omega; H ) }
\nonumber
\\
&&\qquad
= \bigl\| e^{ A {t_1} } \bigl( e^{ A {(t_2 - t_1)} } \xi- \xi\bigr) \bigr\|_{ L^2 ( \Omega;
H ) } \leq\bigl\| e^{ A {(t_2 - t_1 ) } } \xi- \xi\bigr\|_{ L^2 ( \Omega; H ) }
\\
&&\qquad
\leq\bigl\| A^{-1} \bigl( e^{ A {(t_2 - t_1 ) } } - I \bigr) \bigr\|_{ L ( H ) } \| \xi
\|_{ L^2 ( \Omega; D(A) ) } \leq R(t_2 - t_1)\nonumber
\end{eqnarray}
for every $0 \leq t_1 < t_2 \leq T$.
Moreover, we obtain
\begin{eqnarray*}
&& \biggl\| \int_0^{t_2} e^{ A {(t_2 - s )} } F(X_s) \,ds - \int_0^{t_1} e^{ A
{(t_1 - s )} } F(X_s) \,ds \biggr\|_{ L^2 ( \Omega; H ) }
\\
&&\qquad
= \biggl\| \int_{t_1}^{t_2} e^{ A {(t_2 - s )} } F(X_s) \,ds + \int_0^{t_1} \bigl(
e^{ A {(t_2 - s )} } - e^{ A {(t_1 - s )} } \bigr) F(X_s) \,ds \biggr\|_{ L^2 (
\Omega; H ) }
\\
&&\qquad\leq
\int_{t_1}^{t_2} \| F(X_s) \|_{ L^2 ( \Omega; H ) } \,ds + \biggl\|
\int_0^{t_1} \bigl( e^{ A {(t_2 - s )} } - e^{ A {(t_1 - s )} } \bigr) F(X_{s})
\,ds \biggr\|_{ L^2 ( \Omega; H ) }
\end{eqnarray*}
and hence
\begin{eqnarray*}
&& \biggl\| \int_0^{t_2} e^{ A {(t_2 - s )} } F(X_s) \,ds - \int_0^{t_1} e^{ A
{(t_1 - s )} } F(X_s) \,ds \biggr\|_{ L^2 ( \Omega; H ) }
\\
&&\qquad\leq
R (t_2 - t_1) + \int_0^{t_1} \bigl\| e^{ A {(t_2 - s )} } - e^{ A {(t_1 - s
)} } \bigr\|_{ L ( H ) } \| F(X_s) \|_{ L^2 ( \Omega; H ) } \,ds
\\
&&\qquad\leq
R (t_2 - t_1) + R \int_0^{t_1} \bigl\| e^{ A {(t_2 - s )} } - e^{ A {(t_1 -
s )} } \bigr\|^{ (1 - \theta) }_{ L ( H ) } \bigl\| e^{ A {(t_2 - s )} } - e^{ A
{(t_1 - s )} } \bigr\|^\theta_{ L ( H ) } \,ds
\end{eqnarray*}
for every $0 \leq t_1 < t_2 \leq T$.
This yields
\begin{eqnarray*}
&& \biggl\| \int_0^{t_2} e^{ A {(t_2 - s )} } F(X_s) \,ds - \int_0^{t_1} e^{ A
{(t_1 - s )} } F(X_s) \,ds \biggr\|_{ L^2 ( \Omega; H ) }
\\
&&\qquad\leq
R (t_2 - t_1) + 2^{( 1 - \theta)} R \int_0^{t_1} \biggl( \frac{ (t_2 - t_1)
}{ (t_1 - s) } \biggr)^\theta \,ds
\\
&&\qquad\leq
R (t_2 - t_1) + 2 R (t_2 - t_1)^\theta\int_0^{t_1} (t_1 - s)^{-\theta}
\,ds
\end{eqnarray*}
due to Lemma~\ref{lemmaRef4} and therefore, we obtain
%
%
\begin{eqnarray} \label{eqPart2}
&& \biggl\| \int_0^{t_2} e^{ A {(t_2 - s )} } F(X_s) \,ds - \int_0^{t_1} e^{ A
{(t_1 - s )} } F(X_s) \,ds \biggr\|_{ L^2 ( \Omega; H ) } \nonumber
\\
&&\qquad\leq
R (t_2 - t_1) + 2R (t_2 - t_1)^\theta\int_0^{t_1} s^{- \theta} \,ds
\\
&&\qquad\leq
R (t_2 - t_1) + \biggl( \frac{2}{(1-\theta)} R T^{ (1-\theta) } \biggr) (t_2 -
t_1)^\theta\nonumber
\end{eqnarray}
for every $0 \leq t_1 < t_2 \leq T$.
Combining~(\ref{eqPart1}),~(\ref{eqPart2}) and
Assumption~\ref{stochconv} yields the assertion.
\end{pf}
\begin{lemma} \label{lemmaRef3}
Let Assumptions~\ref{semigroup}--\ref{initial} be fulfilled.
Then we obtain
\[
\sup_{0 \leq t_1 < t_2 \leq T} \frac{ \| ( X_{t_2} - O_{t_2} )- (
X_{t_1} - O_{t_1} ) \|_{ L^2 ( \Omega; H ) } } { ( t_2 - t_1 ) } <
\infty,
\]
where
$O\dvtx[0,T] \times\Omega\rightarrow D ( (-A)^\gamma) $
is given in Assumption~\ref{stochconv} and where
$X\dvtx\break[0,T] \times\Omega\rightarrow D ( (-A)^\gamma) $
is the solution of the SPDE~(\ref{eqSolution}).
\end{lemma}
\begin{pf}
First, let $R \in[0, \infty)$ be the real number given by
\[
R := \| \xi\|_{ L^2 ( \Omega; D(A) ) } + \sup_{ t \in[0,T] } \|
F(X_t) \|_{ L^2 ( \Omega; H ) } + \sup_{0 \leq t_1 < t_2 \leq T} \frac{
\| X_{t_2} - X_{t_1} \|_{ L^2 ( \Omega; H ) } } { ( t_2 - t_1 )^\theta
},
\]
which exists due to Lemma~\ref{lemmaRef2}.
Then we have
%
%
\begin{eqnarray} \label{eqRef2}
\| e^{ A {t_2} } \xi- e^{ A {t_1} } \xi\|_{ L^2 ( \Omega; H ) } &=& \bigl\|
e^{ A {t_1} } \bigl( e^{ A {(t_2 - t_1)} } \xi- \xi\bigr) \bigr\|_{ L^2 ( \Omega; H )}
\nonumber\\
&\leq&\bigl\| e^{ A {(t_2 - t_1 ) } } \xi- \xi\bigr\|_{ L^2 ( \Omega; H ) } \\
&\leq&
R (t_2 - t_1)\nonumber
\end{eqnarray}
for every $0 \leq t_1 < t_2 \leq T$. %
Moreover, we have
\begin{eqnarray*}
&& \biggl\| \int_0^{t_2} e^{ A {(t_2 - s )} } F(X_s) \,ds - \int_0^{t_1} e^{ A
{(t_1 - s )} } F(X_s) \,ds \biggr\|_{ L^2 ( \Omega; H ) }
\\
&&\qquad
= \biggl\| \int_{t_1}^{t_2} e^{ A {(t_2 - s )} } F(X_s) \,ds\\
&&\qquad\quad\hspace*{6pt} + \int_0^{t_1} \bigl(
e^{ A {(t_2 - s )} } - e^{ A {(t_1 - s )} } \bigr) F(X_s) \,ds \biggr\|_{ L^2 (
\Omega; H ) }
\\
&&\qquad\leq
\int_{t_1}^{t_2} \bigl\| e^{A (t_2 - s )} \bigr\|_{ L ( H ) } \| F(X_s) \|_{ L^2
( \Omega; H ) } \,ds
\\
&&\qquad\quad{}
+ \biggl\| \int_{0}^{t_1} \bigl( e^{A (t_2 - s )} -e^{A (t_1 - s )} \bigr) \bigl(F(X_s) -
F(X_{t_1})\bigr) \,ds \biggr\|_{ L^2 ( \Omega; H ) }
\\
&&\qquad\quad{}
+ \biggl\| \int_0^{t_1} \bigl( e^{A (t_2 - s )} - e^{A (t_1 - s )} \bigr) F(X_{t_1}) \,ds
\biggr\|_{ L^2 ( \Omega; H ) }
\end{eqnarray*}
and therefore
\begin{eqnarray*}
&& \biggl\| \int_0^{t_2} e^{A (t_2 - s )} F(X_s) \,ds - \int_0^{t_1} e^{A (t_1 -
s )} F(X_s) \,ds \biggr\|_{ L^2 ( \Omega; H ) }
\\
&&\qquad\leq
R (t_2 - t_1)\\
&&\qquad\quad{} + \int_{0}^{t_1} \bigl\| e^{A (t_2 - s )} -e^{A (t_1 - s )}
\bigr\|_{ L ( H ) } \| F(X_s) - F(X_{t_1}) \|_{ L^2 ( \Omega; H ) } \,ds
\\
&&\qquad\quad{}
+ R \biggl\| \int_0^{t_1} \bigl( e^{A (t_2 - s )} - e^{A (t_1 - s )} \bigr) \,ds \biggr\|_{ L (
H ) }
\end{eqnarray*}
for every $0 \leq t_1 < t_2 \leq T$.
Hence, we obtain
\begin{eqnarray*}
&& \biggl\| \int_0^{t_2} e^{A (t_2 - s )} F(X_s) \,ds - \int_0^{t_1} e^{A (t_1 -
s )} F(X_s) \,ds \biggr\|_{ L^2 ( \Omega; H ) }
\\
&&\qquad\leq
R (t_2 - t_1)\\
&&\qquad\quad{} + c \int_{0}^{t_1} \bigl\| e^{A (t_2 - s )} -e^{A (t_1 - s )}
\bigr\|_{ L ( H ) } \| X_s - X_{t_1} \|_{ L^2 ( \Omega; H ) } \,ds
\\
&&\qquad\quad{}
+ R \biggl\| \int_0^{t_1} e^{ A {((t_2 - t_1) + s )} } \,ds - \int_0^{t_1} e^{A
s} \,ds \biggr\|_{ L ( H ) }
\end{eqnarray*}
and
\begin{eqnarray*}
&& \biggl\| \int_0^{t_2} e^{A (t_2 - s )} F(X_s) \,ds - \int_0^{t_1} e^{A (t_1 -
s )} F(X_s) \,ds \biggr\|_{ L^2 ( \Omega; H ) }
\\
&&\qquad\leq
R (t_2 - t_1) + cR \int_{0}^{t_1} \bigl\| e^{A (t_2 - s )} -e^{A (t_1 - s )}
\bigr\|_{ L ( H ) } | s - t_1 |^\theta \,ds
\\
&&\qquad\quad{}
+ R \biggl\| \int_{(t_2 - t_1)}^{t_2} e^{A s} \,ds - \int_0^{t_1} e^{A s} \,ds
\biggr\|_{ L ( H ) }
\end{eqnarray*}
for every $0 \leq t_1 < t_2 \leq T$.
This shows
\begin{eqnarray*}
&&\biggl\| \int_0^{t_2} e^{A (t_2 - s )} F(X_s) \,ds - \int_0^{t_1} e^{A (t_1 -
s )} F(X_s) \,ds \biggr\|_{ L^2 ( \Omega; H ) } \\
&&\qquad\leq R (t_2 - t_1)
\\
&&\qquad\quad{}
+ cR \int_{0}^{t_1} \frac{(t_2 - t_1)}{(t_1 - s)} | s - t_1 |^\theta \,ds
+ R \biggl\| \int_{t_1}^{t_2} e^{A s} \,ds - \int_0^{(t_2 - t_1)} e^{A s} \,ds
\biggr\|_{ L ( H ) }
\end{eqnarray*}
due to Lemma~\ref{lemmaRef4} and
\begin{eqnarray*}
&& \biggl\| \int_0^{t_2} e^{A (t_2 - s )} F(X_s) \,ds - \int_0^{t_1} e^{A (t_1 -
s )} F(X_s) \,ds \biggr\|_{ L^2 ( \Omega; H ) }
\\
&&\qquad
\leq R (t_2 - t_1) + cR (t_2 - t_1) \int_{0}^{t_1} (t_1 - s)^{(\theta-
1)} \,ds + 2 R (t_2 - t_1)
\\
&&\qquad
= R (t_2 - t_1) + cR (t_2 - t_1) \int_{0}^{t_1} s^{(\theta- 1)} \,ds + 2
R (t_2 - t_1)
\\
&&\qquad
\leq\bigl(R + cR (T+1) \theta^{-1} + 2 R\bigr) (t_2 - t_1)
\end{eqnarray*}
for every $0 \leq t_1 < t_2 \leq T$.
Combining this and~(\ref{eqRef2}) shows the assertion.
\end{pf}

\section*{Acknowledgments}
We strongly thank the anonymous referee for his careful
reading and his very valuable advice.

\printaddresses

\end{document}